\documentclass[a4paper,12pt]{article}
\usepackage{amssymb}
\usepackage{amsfonts}
\usepackage{amsmath}
\usepackage{color}
\usepackage{latexsym,array,theorem,bbm}
\usepackage[colorlinks=true,linkcolor=blue,citecolor=red]{hyperref}    

\sffamily
\newcommand{\RR}{\mathbb{R}}
\newcommand{\ZZ}{\mathbb{Z}}

\newcommand{\NN}{\mathbb{N}}

\newcommand{\cA}{{\mathcal A}}
\newcommand{\cD}{{\mathcal D}}
\newcommand{\cL}{{\mathcal L}}

\renewcommand{\P}{\mathsf{P}}
\newcommand{\cov}{\mathsf{Cov}}
\newcommand{\EE}{\mathsf{E}}
\newcommand{\PP}{\mathsf{P}}
\newcommand{\dd}{\mathsf{d}}  
\newcommand{\ee}{\mathrm{e}}  
\newcommand{\Prob}{\mathsf{P}}

\newcommand{\aWb}{$\alpha$\hspace{-1.5pt}-\hspace{-0.5pt}Wiener bridge}
\newcommand{\proofend}{\hfill\mbox{$\Box$}\medskip}

\newcommand{\distre}{\stackrel{\cL}{=}}

\numberwithin{equation}{section}

\theoremstyle{change}
\newtheorem{Lem}{Lemma.}[section]
\newtheorem{Thm}[Lem]{Theorem.}
\newtheorem{Pro}[Lem]{Proposition.}
\newtheorem{Cor}[Lem]{Corollary.}
\theorembodyfont{\rm}

\newcounter{Not}

\newtheorem{Rem}[Lem]{Remark.}

\begin{document}

\begin{center}
 {\bfseries\Large Karhunen--Lo\`eve expansions of \aWb s}\\[5mm]

 {\sc\large M\'aty\'as $\text{Barczy}^{*,\diamond}$} {\large and}
 {\sc\large Endre $\text{Igl\'oi}^*$}
\end{center}

\vskip0.2cm

* University of Debrecen, Faculty of Informatics, Pf.~12, H--4010 Debrecen, Hungary;
 e--mails: barczy.matyas@inf.unideb.hu (M. Barczy), igloi@tigris.unideb.hu (E. Igl\'oi).

$\diamond$ Corresponding author.

\vskip0.2cm

\renewcommand{\thefootnote}{}

\vspace*{-5mm}

\date{}

\begin{abstract}
We study Karhunen--Lo\`eve expansions of the process $(X_t^{(\alpha)})_{t\in[0,T)}$ given by the
 stochastic differential equation
 $\dd X_t^{(\alpha)} = -\frac\alpha{T-t}\,X_t^{(\alpha)}\dd t+\dd B_t,$  $t\in[0,T),$ with an initial condition
 $X_0^{(\alpha)}=0,$ where $\alpha>0,$ $T\in(0,\infty)$ and $(B_t)_{t\geq 0}$ is a standard Wiener process.
This process is called an $\alpha$-Wiener bridge or a scaled Brownian bridge, and
 in the special case of $\alpha=1$ the usual Wiener bridge.
We present weighted and unweighted Karhunen--Lo\`eve expansions of $X^{(\alpha)}$.
As applications, we calculate the Laplace transform and the distribution function of
 the $L^2[0,T]$\!-norm square of $X^{(\alpha)}$ \ studying also its asymptotic behavior
 (large and small deviation).
\end{abstract}

\section{Introduction}
\label{Isect}

There are few stochastic processes of interest, even among Gaussian ones, for which the Karhunen--Lo\`eve (KL) expansion is explicitly known.
Some examples are those of the Wiener process, the Ornstein--Uhlenbeck process and the Wiener bridge, see, e.g.,
 Ash and Gardner \cite[Example 1.4.4]{AshGar}, Papoulis \cite[Problem 12.7]{Papo},
 Liu and Ulukus \cite[Section III]{LiuUlu}, Corlay and Pag\`{e}s \cite[Section 5.4 B]{CorPag}
 and Deheuvels \cite[Remark 1.1]{Deh1}.
Recently, there is a renewed interest  in this field:
 some KL expansions were provided for weighted Wiener processes and weighted Wiener bridges
 with weighting function having the form $t^\beta$ (these expansions make use of Bessel functions),
 see Deheuvels and Martynov \cite{DehMar}.
The most recent results on this topic are those of Deheuvels, Peccati and Yor \cite{DehPecYor},
 Deheuvels \cite{Deh1}, \cite{Deh2}, Luschgy and Pag\`es \cite{LusPag},
 Nazarov and Nikitin \cite{NazNik} and Nazarov and Pusev \cite{NazPus} (the latter two ones are about
 exact small deviation asymptotics for weighted $L^2$-norm of some Gaussian processes).

Let $0<S<T<\infty$ and $0<\alpha<\infty$ be arbitrarily fixed
and let $(B_t)_{t\geq0}$ be a standard Wiener process on a probability space $(\Omega,\cal{A},\P).$
Let us consider the stochastic differential equation (SDE)
\[
\left\{
\begin{aligned}
 \dd X_t^{(\alpha)}&=-\frac\alpha{T-t}\,X_t^{(\alpha)}\dd t+\dd B_t,\qquad t\in[0,S],\\
 X_0^{(\alpha)}&=0.
\end{aligned}
\right.
\]
The drift and diffusion terms satisfy the usual Lipschitz and linear growth conditions,
 so, by {\O}ksendal \cite[Theorem 5.2.1]{Oks} or Jacod and Shiryaev \cite[Chapter III, Theorem 2.32]{JacShi},
 the above SDE has a strong solution which is pathwise unique (i.e., it has a unique strong solution).
Since $S\in(0,T)$ is chosen arbitrarily, we obtain by successive extension that also the SDE
\begin{align}
\label{alpha_W_bridge}
\left\{
\begin{aligned}
 \dd X_t^{(\alpha)}&=-\frac\alpha{T-t}\,X_t^{(\alpha)}\dd t+\dd B_t,\qquad t\in[0,T),\\
 X_0^{(\alpha)}&=0,
\end{aligned}
\right.
\end{align}
has a unique strong solution.
Namely, it is
\begin{equation}
\label{MArepr}
X_t^{(\alpha)} = \int_0^t\!\bigg(\frac{T-t}{T-s}\bigg)^\alpha \dd B_s,\qquad t\in[0,T),
\end{equation}
as it can be checked by It\^o's formula.
The process $(X^{(\alpha)}_t)_{t\in[0,T)}$ given by \eqref{MArepr} is called an
 \aWb\ (from $0$ to $0$ on the time interval $[0,T]$\!).
To our knowledge, these kind of processes have  been first considered by
 Brennan and Schwartz \cite{BreSch}, and see also Mansuy \cite{Man}.
In Brennan and Schwartz \cite{BreSch} \aWb s are used to model
 the arbitrage profit associated with a given futures contract in the absence of transaction costs.
Sondermann, Trede and Wilfling  \cite{SonTreWil} and Trede and Wilfling \cite{TreWil}
use the SDE (\ref{alpha_W_bridge}) to describe the fundamental component of an exchange rate process
 and they call the process $X^{(\alpha)}$ as a scaled Brownian bridge.
The essence of these models is that
 the coefficient of $X^{(\alpha)}_t$ in the drift term in (\ref{alpha_W_bridge}) represents some kind of mean
 reversion, a stabilizing force that keeps pulling the process towards its mean (zero in this reduced form),
 and the absolute value of this force is  increasing proportionally to the inverse of the remaining time $T-t,$
 with the rate constant $\alpha$.

This process has been also studied by Barczy and Pap \cite{BarPap1}, \cite[Section 4]{BarPap2} from several
 points of view, e.g., singularity of probability measures induced by the process $X^{(\alpha)}$
 with different values of $\alpha,$ sample path properties, Laplace transforms of certain functionals
 of $X^{(\alpha)}$ and maximum likelihood estimation of $\alpha$.
The process $(X^{(\alpha)}_t)_{t\in[0,T)}$ is  Gaussian and
 for all $t\in[0,T),$ $\EE X^{(\alpha)}_t=0$ and the covariance function of $X^{(\alpha)}$
 given in Barczy and Pap \cite[Lemma 2.1]{BarPap1} is
\begin{align}\label{covfv}
 \begin{split}
 R^{(\alpha)}(s,t)&:=\cov\big(X_s^{(\alpha)},X_t^{(\alpha)}\big) \\[1mm]
  &=\begin{cases}\frac{(T-s)^\alpha(T-t)^\alpha}{1-2\alpha}\big(T^{1-2\alpha}-(T-(s\wedge t))^{1-2\alpha}\big)
  &\text{if }\alpha\neq\frac12,\\[6pt]
    \sqrt{(T-s)(T-t)}\,\ln\!\big(\frac{T}{T-(s\wedge t)}\big)
  &\text{if }\alpha=\frac12,
 \end{cases}
 \end{split}
\end{align}
 for all $s,t\in[0,T),$ where $s\wedge t:=\min(s,t)$.
By Barczy and Pap \cite[Lemma 3.1]{BarPap1}, the \aWb $(X_t^{(\alpha)})_{t\in[0,T)}$
 has an almost surely continuous extension $(X^{(\alpha)}_t)_{t\in[0,T]}$ to the time interval $[0,T]$ such that
 $X_T^{(\alpha)}=0$ with probability one.
The possibility of such an extension is based on that the parameter $\alpha$ is positive and
 on a strong law of large numbers for square integrable local martingales.
We note here also that (\ref{alpha_W_bridge}--\ref{MArepr}) continue to hold for $\alpha\leq0$ as well.
However, there does not exist an almost surely continuous extension of the process
 $(X^{(\alpha)}_t)_{t\in[0,T)}$ onto $[0,T]$ which would take some constant at time $T$
 with probability one (i.e., which would be a bridge), and this is why the range of
 the parameter $\alpha$ is restricted to positive values.
Indeed, for $\alpha=0$ we obtain the Wiener process,
and in case of $\alpha<0$ the second moment of the solution $X^{(\alpha)}_t$ given by (\ref{MArepr}) converges to infinity, as (\ref{covfv}) (with $s=t$\!) shows. Hence the assumption of the existence of an almost surely continuous extension to $[0,T]$ such that this extension takes some constant at time $T$ with probability one
 (i.e., we have a bridge) would result in a contradiction.
We note that another proof of the impossibility of such an extension in the case of $\alpha<0$ can be found
 in Barczy and Pap \cite[Remark 3.5]{BarPap1}.
Finally, we remark that Mansuy \cite[Proposition 4]{Man} studied the question whether it is possible
 to derive the \aWb\ from a (single) Gaussian process by taking a bridge.

Next we check that the \aWb $(X_t^{(\alpha)})_{t\in[0,T]}$ is $L^2$-continuous.
By Theorem 1.3.4 in Ash and Gardner \cite{AshGar}, it is enough to show that the covariance
 function $R^{(\alpha)}$ can be extended continuously onto $[0,T]^2:=[0,T]\times[0,T]$
 such that this extension (which will be also denoted by $R^{(\alpha)}$) is zero on the set
 $\{(s,T) : s\in[0,T]\}\cup\{(T,t) : t\in[0,T]\}$.
This follows by
 \begin{equation}
 \label{vani}
   \lim_{(s,t)\to(s_0,T)}R^{(\alpha)}(s,t)=\lim_{(s,t)\to(T,t_0)}R^{(\alpha)}(s,t)=0,
     \qquad s_0,t_0\in[0,T].
 \end{equation}
Indeed, if $\alpha\ne1/2$ and $s_0<T,$ then
 \begin{align*}
  \lim_{(s,t)\to(s_0,T)}R^{(\alpha)}(s,t)
    = \frac{(T-s_0)^\alpha}{1-2\alpha}(T^{1-2\alpha} - (T-s_0)^{1-2\alpha})
      \lim_{t\uparrow T} (T-t)^\alpha
    =0.
 \end{align*}
If $0<\alpha<1/2$ and $s_0=T,$ then
 \begin{align*}
  \lim_{(s,t)\to(s_0,T)}R^{(\alpha)}(s,t)
    = \frac{T^{1-2\alpha}}{1-2\alpha}
      \lim_{(s,t)\uparrow(T,T)} (T-s)^\alpha (T-t)^\alpha
    =0.
 \end{align*}
If $\alpha>1/2$ and $s_0=T,$ then
 \begin{align*}
  \lim_{(s,t)\to(s_0,T)}R^{(\alpha)}(s,t)
   & = \frac{1}{2\alpha-1}
      \lim_{(s,t)\uparrow(T,T)} \frac{(T-s)^\alpha (T-t)^\alpha}{(T- (s\wedge t))^{2\alpha-1}}\\
   & \leq \frac{1}{2\alpha-1}
         \lim_{(s,t)\uparrow(T,T)} (T-(s\wedge t))
     =0.
 \end{align*}
If $\alpha=1/2$ and $s_0<T,$ then
 \begin{align*}
  \lim_{(s,t)\to(s_0,T)}R^{(\alpha)}(s,t)
    = \sqrt{T-s_0}\ln\left(\frac{T}{T-s_0}\right)
      \lim_{t\uparrow T} \sqrt{T-t}
    =0.
 \end{align*}
If $\alpha=1/2$ and $s_0=T,$ then
 \begin{align*}
  \lim_{(s,t)\to(s_0,T)}R^{(\alpha)}(s,t)
    & = \lim_{(s,t)\uparrow(T,T)}
        \sqrt{(T-s)(T-t)}\ln\left(\frac{T}{T-(s\wedge t)}\right)\\
    & \leq \lim_{(s,t)\uparrow(T,T)}
         (T-(s\wedge t))\ln\left(\frac{T}{T-(s\wedge t)}\right)
    =0.
 \end{align*}

We also have $R^{(\alpha)}\in L^2([0,T]^2)$.
So, the integral operator associated to the kernel function $R^{(\alpha)},$ i.e., the operator
 $A_{R^{(\alpha)}}:L^2[0,T]\rightarrow L^2[0,T],$
\begin{equation}
\label{Rdefi}
\begin{aligned}
 (A_{R^{(\alpha)}}(\phi))(t)&:=\int_0^T\! R^{(\alpha)}(t,s)\phi(s)\,\dd s,\quad t\in[0,T],
  \qquad \phi\in L^2[0,T],
\end{aligned}
\end{equation}
is of the Hilbert--Schmidt type, thus $(X_t^{(\alpha)})_{t\in[0,T]}$ has a Karhunen--Lo\`eve (KL) expansion
based on $[0,T]$\!:
\begin{equation}
\label{KLrepr}
 X_t^{(\alpha)} =\sum_{k=1}^\infty
   \sqrt{\lambda^{(\alpha)}_k}\,\xi_k e^{(\alpha)}_k(t),\quad t\in[0,T],
\end{equation}
 where $\xi_k,\ k\in\mathbb{N},$ are independent standard normally distributed random variables,
$\lambda^{(\alpha)}_k,\ k\in\mathbb{N},$ are the non-zero eigenvalues of the integral
 operator $A_{R^{(\alpha)}}$ and $ e^{(\alpha)}_k(t),\ t\in[0,T],\ k\in\mathbb{N},$ are the corresponding normed eigenfunctions, which are pairwise orthogonal in $L^2[0,T],$
see, e.g., Ash and Gardner \cite[Theorem 1.4.1]{AshGar}.
Observe that (\ref{KLrepr}) has infinitely many terms. Indeed, if it had a finite number of terms,
i.e., if there were only a finite number of eigenfunctions, say $N,$ then by the help of (\ref{alpha_W_bridge})
 (considering it as an integral equation) we would obtain that the Wiener process $(B_t)_{t\in[0,T]}$ is
 concentrated in an $N$ dimensional subspace of $L^2[0,T],$ and so even of $C[0,T],$ with probability one.
This results in a contradiction, since the integral operator associated to the covariance function
 (as a kernel function) of a standard Wiener process has infinitely many eigenvalues and eigenfunctions.
We also note that the normed eigenfunctions
are unique only up to sign.
The series in (\ref{KLrepr}) converges in $L^2(\Omega,\cA,P)$ to $X_t^{(\alpha)},$
 uniformly on $[0,T],$ i.e.,
 \[
    \sup_{t\in[0,T]}\EE\left( \left\vert X_t^{(\alpha)}
     - \sum_{k=1}^n\sqrt{\lambda^{(\alpha)}_k}\,\xi_k e^{(\alpha)}_k(t)\right\vert^2\right)\to 0
     \qquad \text{as $n\to\infty.$}
 \]
Moreover, since $R^{(\alpha)}$ is continuous on $[0,T]^2,$  the eigenfunctions corresponding
 to non-zero eigenvalues are also continuous on $[0,T],$ see, e.g. Ash and Gardner \cite[p. 38]{AshGar}
 (this will be important in the proof of Theorem \ref{KLthm}, too).
Since the terms on the right-hand side of (\ref{KLrepr}) are independent normally distributed random variables
 and $(X_t^{(\alpha)})_{t\in[0,T]}$ has continuous sample paths with probability one, the series converges
 even uniformly on $[0,T]$ with probability one (see, e.g., Adler \cite[Theorem 3.8]{Adler}).

The rest of the paper is organized as follows.
In Section \ref{KLsect} we make the KL representation \eqref{KLrepr} of the \aWb\ $(X_t^{(\alpha)})_{t\in[0,T]}$
 as explicit as it is possible, see Theorem \ref{KLthm}.
We also consider two special cases of the KL representation \eqref{KLrepr}, first we study the case
 $\alpha\downarrow 0$ (standard Wiener process) and then the case $\alpha=1$ (Wiener bridge),
 see Remark \ref{KL0} and Remark \ref{KL1}, respectively.
In Remark \ref{KL2} we present a so-called space-time transformed Wiener process representation of the \aWb.
Using this representation and a result of Guti\'errez and Valderrama
 \cite[Theorem 1]{GutVal}, we derive a weighted KL representation of the \aWb
 $(X_t^{(\alpha)})_{t\in[0,T]}$ based on $[0,S],$ where $0<S<T,$ see Theorem \ref{wghtKLthm}.
We also consider two special cases of this weighted KL representation, first we study the case $\alpha\downarrow 0$ and then the case $\alpha=1,$ see Remark \ref{wghtKLrem1} and
 Remark \ref{wghtKLrem2}, respectively.
Further, we give an infinite series representation of $\int_0^S(X^{(1/2)}_u)^2/(T-u)^2\,\dd u,$
 where $0<S<T$, see Remark \ref{wghtKLrem3}.
Section \ref{appl} is devoted to the applications.
In Proposition \ref{appl1} we determine the Laplace transform of
 the $L^2[0,T]$\!-norm square of $(X_t^{(\alpha)})_{t\in[0,T]}$
 and of the $L^2[0,S]$\!-norm square of $(X_t^{(1/2)}/(T-t))_{t\in[0,S]},$ where $0<S<T.$
In Corollary \ref{COR:Bessel} we give a new probabilistic proof for the well-known result of
 Rayleigh, namely, for the sum of the square of the reciprocals of the positive zeros of Bessel
 functions of the first kind (with order greater than $-1/2$ in our case).
Based on the Smirnov formula (see, e.g., Smirnov \cite[formula (97)]{Smi}) we write the survival function
 of the $L^2[0,T]$\!-norm square of $(X_t^{(\alpha)})_{t\in[0,T]}$ in an infinite series form,
 see Proposition \ref{appl2}.
We also consider two special cases of Proposition \ref{appl2}, the case $\alpha\downarrow 0$
 and the case $\alpha=1,$ see Remark \ref{appl3} and Remark \ref{appl4}.
In Corollary \ref{appl5}, based on a result of Zolotarev \cite{Zol}, we study large deviation
 probabilities for the $L^2[0,T]$\!-norm square of the $\alpha$-Wiener bridge.
Finally, based on a result of Li \cite[Theorem 2]{Li1}, we describe the behavior at zero
 of the distribution function (small deviation probabilities) of the $L^2[0,T]$\!-norm square of the
 $\alpha$-Wiener bridge, see Corollay \ref{appl6}.
In the appendix we list some important properties of Bessel functions of the first and second kind,
 respectively.

We remark that our results for \aWb s may have some generalizations
for random fields.
Namely, for all $S>0,$ $T>0$ and $\alpha\geq 0,$ $\beta\geq 0,$ one can consider a zero-mean
 Gaussian random field $(X_{s,t}^{(\alpha,\beta)})_{(s,t)\in[0,S]\times[0,T]}$ with the covariances
 \[
   \EE(X_{s_1,t_1}^{(\alpha,\beta)}X_{s_2,t_2}^{(\alpha,\beta)})
       = R^{(\alpha)}(s_1,s_2) R^{(\beta)}(t_1,t_2),
         \qquad (s_1,t_1), (s_2,t_2)\in [0,S]\times[0,T].
 \]
Such a random field exists, since $(X_s^{(\alpha)}X_t^{(\beta)})_{(s,t)\in[0,S]\times[0,T]}$ admits the above
 covariances, where $X^{(\alpha)}$ and $X^{(\beta)}$ are independent, and Kolmogorov's consistency theorem
 comes into play.
This class of Gaussian processes may deserve more attention since it would generalize some well-known
 limit processes in mathematical statistics such as the Kiefer process (known also a tied down Brownian
 sheet), see, e.g., Cs\"org\H o and R\'ev\'esz \cite[Section 1.15]{CsorRev}.
Indeed, with $S=1,$ $T=\infty$, $\alpha=1$ and $\beta=0$ the process $X^{(\alpha,\beta)}$
 is nothing else but the Kiefer process having covariance function $(s_1\wedge s_2 - s_1s_2)(t_1\wedge t_2),$
 $(s_1,t_1),(s_2,t_2)\in[0,1]\times[0,\infty).$

In all what follows $\NN$, $\ZZ_+$ and $\ZZ$ denote the set of natural numbers, nonnegative
 integers and integers, respectively.

\section{Karhunen--Lo\`eve expansions of $\alpha$-Wiener bridges}
\label{KLsect}

First we recall the notion of Bessel functions of the first kind which plays a key role
 in the KL expansions we will obtain, and also the Bessel functions of the second kind
 (also called Neumann functions) which will appear in the proofs.
They can be defined, resp., as
 \begin{align*}
 &J_\nu(x):=\sum_{k=0}^\infty\frac{(-1)^k}{k!\,\Gamma(k+\nu+1)}\,\left(\frac{x}2\right)^{2k+\nu},
             \quad x\in(0,\infty),\;\nu\in\RR,\\[5pt]
 &Y_\nu(x):=
         \begin{cases}
            \frac{J_\nu(x)\cos(\pi\nu)-J_{-\nu}(x)}{\sin(\pi\nu)},
               \quad & \text{if $x\in(0,\infty),\;\nu\in\RR\setminus\ZZ$},\\[1mm]
            Y_n(x):=\lim_{\mu\to n}Y_\mu(x)
               \quad & \text{if $x\in(0,\infty),$ $\nu=n\in\ZZ,$}
          \end{cases}
 \end{align*}
 where $\Gamma(z)$ for $z<0, \, z\not\in\ZZ,$ is defined by a recursive application of the rule
 $\Gamma(z)=\Gamma(z+1)/z$, $z<0$, $z\not\in\ZZ$, and we use the convention that $1/\Gamma(-k):=0$, $k\in\ZZ_+,$
 yielding that the first $n$ terms in the series of $J_\nu(x)$ vanish if $\nu=-n,\,n\in\NN$,
 see, e.g., Watson \cite[pp. 40, 64]{Wat}.

In all what follows we will put $\nu:=\alpha-1/2,$ where $\alpha>0$.
Next we present our main theorem.

\begin{Thm}
\label{KLthm}
Let $\alpha>0,$ $\nu:=\alpha-1/2,$ and $z^{(\nu)}_k\!,\ k\in\mathbb{N},$ be the (positive) zeros
 of $J_\nu.$
Then in the KL expansion (\ref{KLrepr}) of the \aWb\ $(X^{(\alpha)}_t)_{t\in[0,T]}$
the eigenvalues and the corresponding normed eigenfunctions are
\begin{align}
 &\qquad\qquad\lambda^{(\alpha)}_k = \frac{T^2}{(z_k^{(\nu)})^2},\quad  k\in\mathbb{N},\label{lambdakk}\\[2mm]
 e^{(\alpha)}_k(t) &=\sqrt{\frac{2}T\left(1-\frac{t}T\right)}\,
 \frac{J_\nu\big(z_k^{(\nu)}(1-t/T)\big)}{\big\vert J_{\nu+1}\big(z_k^{(\nu)}\big)\big\vert}\notag\\[2mm]
 &=\sqrt{\frac{2}T\left(1-\frac{t}T\right)}\,
 \frac{J_\nu\big(z_k^{(\nu)}(1-t/T)\big)}{\big\vert J_{\nu-1}\big(z_k^{(\nu)}\big)\big\vert},
 \quad t\in[0,T],\notag
\end{align}
where we take the continuous extension of $e_k^{(\alpha)}$ at $t=T$ for $-1/2<\nu<0$,
 i.e., $e_k^{(\alpha)}(T)=0$ for $\alpha<1/2$ (see part (i) of Proposition \ref{PRO:Bessel}).
\end{Thm}

\noindent\textbf{Proof.}
For simplicity we first assume that $T=1$. We will return to the interval $[0,T]$ at the end of the proof.

We consider the cases $\alpha\ne1/2$ and $\alpha=1/2$ separately, since the covariance function
 $R^{(\alpha)}$ has different forms in these two cases.

If  $\alpha\neq1/2,$ i.e., $\nu\neq0,$  then the covariance function of $(X^{(\alpha)}_t)_{t\in[0,1]}$
 takes the form
\begin{equation}
\label{Ralpha}
R^{(\alpha)}(t,s)
 =\begin{cases}
 \frac{(1-s)^\alpha(1-t)^\alpha}{1-2\alpha}\left(1-(1-(s\wedge t))^{1-2\alpha}\right)
   &\text{if \ $s,t\in[0,1)$},\\[2pt]
  0&\text{if \ $s\in\{0,1\}$ \ or \ $t\in\{0,1\}$}.
\end{cases}
\end{equation}
Let $\lambda$ be a non-zero (and hence positive) eigenvalue of the integral operator $A_{R^{(\alpha)}}$.
Then the eigenvalue equation is
\begin{equation}\label{eigeq}
\int_0^1\! R^{(\alpha)}(t,s)e(s)\,\dd s=\lambda e(t),\quad t\in[0,1].
\end{equation}
With the special choices $t=0$ and $t=1$ we have the boundary conditions
 $e(0)= e(1)=0$.
By \eqref{eigeq}, we have
 \[
   \int_0^t R^{(\alpha)}(t,s)e(s)\,\dd s + \int_t^1 R^{(\alpha)}(t,s)e(s)\,\dd s
                      =\lambda e(t),\quad t\in(0,1),
 \]
 and differentiating twice with respect to $t$ (which can be done since the integrands
 are continuously differentiable with respect to $t$) gives
\begin{align}\label{DE_e1}
\lambda e'(t)&=-\frac{\alpha\lambda}{1-t}\, e(t)
       + (1-t)^{-\alpha}\int_t^1\!(1-s)^\alpha e(s)\,\dd s,\\[5pt] \label{DE_e2}
\lambda e''(t)&=-\frac{\alpha\lambda}{(1-t)^2}\, e(t)-\frac{\alpha\lambda}{1-t}\, e'(t)+\alpha(1-t)^{-1-\alpha}\int_t^1\!(1-s)^\alpha e(s)\,\dd s- e(t)
\end{align}
 for all $t\in(0,1)$.
Indeed, by \eqref{eigeq}, for all $t\in(0,1),$
 \[
   \int_0^t\! R^{(\alpha)}(t,s)e(s)\,\dd s + \int_t^1\! R^{(\alpha)}(t,s)e(s)\,\dd s = \lambda e(t),
 \]
 and then
 \begin{align}\label{seged2}
  \begin{split}
   \lambda e(t)
    & = \frac{(1-t)^\alpha}{1-2\alpha} \int_0^t \big((1-s)^\alpha -(1-s)^{1-\alpha}\big)e(s)\,\dd s\\[1mm]
    & \phantom{=\;}
       + \frac{(1-t)^\alpha - (1-t)^{1-\alpha}}{1-2\alpha} \int_t^1 (1-s)^\alpha e(s)\,\dd s,
       \quad t\in(0,1),
   \end{split}
 \end{align}
 or equivalently
 \begin{align}\label{seged1}
  \begin{split}
 \frac{(1-2\alpha)}{(1-t)^\alpha} \lambda e(t)
  & = \int_0^1 (1-s)^\alpha e(s)\,\dd s
           - \int_0^t (1-s)^{1-\alpha} e(s)\,\dd s \\
  & \phantom{=\;}
           - (1-t)^{1-2\alpha}\int_t^1 (1-s)^\alpha e(s)\,\dd s,
  \quad t\in(0,1).
  \end{split}
 \end{align}
Differentiating \eqref{seged2} with respect to $t$ we have
 \begin{align*}
  \lambda e'(t)
    & = -\frac{\alpha (1-t)^{\alpha-1}}{1-2\alpha}
      \int_0^t \big((1-s)^\alpha - (1-s)^{1-\alpha}\big) e(s)\,\dd s \\[1mm]
    & \phantom{=\;} + \frac{-\alpha(1-t)^{\alpha-1} + (1-\alpha)(1-t)^{-\alpha}}{1-2\alpha}
                \int_t^1 (1-s)^\alpha e(s)\,\dd s ,
           \quad t\in(0,1),
 \end{align*}
 or equivalently
 \begin{align*}
   \lambda e'(t)
    &= -\frac{\alpha (1-t)^{\alpha-1}}{1-2\alpha} \int_0^1 (1-s)^\alpha e(s)\,\dd s
         + \frac{\alpha (1-t)^{\alpha-1}}{1-2\alpha} \int_0^t (1-s)^{1-\alpha} e(s)\,\dd s \\[1mm]
    & \phantom{=\;}
         + \frac{(1-\alpha)(1-t)^{-\alpha}}{1-2\alpha}\int_t^1 (1-s)^{\alpha} e(s)\,\dd s ,
         \quad t\in(0,1).
 \end{align*}
Using \eqref{seged1}, one can derive \eqref{DE_e1}.
Combining the equations \eqref{DE_e1} and \eqref{DE_e2} we have
\begin{equation}
\label{RBdiffeq}
 e''(t)=-\left(\frac1\lambda+\frac{\alpha(1-\alpha)}{(1-t)^2}\right) e(t),
  \quad t\in(0,1).
\end{equation}
Let us introduce the function $f:[0,1/\sqrt\lambda]\to\RR,$
 \[
     f(x):= e(1-\sqrt\lambda x),\qquad x\in[0,1/\sqrt\lambda].
 \]
By \eqref{RBdiffeq}, we obtain the Riccati--Bessel differential equation
 (first studied by Plana, see Watson \cite[p. 95]{Wat})
\begin{equation}
\label{RBdifeq}
x^2 f''(x)+\left(x^2-(\alpha-1)\alpha\right) f(x)=0,\qquad x\in(0,1/\sqrt\lambda)
\end{equation}
with the boundary conditions $ f(0)= f(1/\!\sqrt\lambda\,)=0.$

In the other case, when  $\alpha=1/2,$ i.e., $\nu=0,$
the covariance function of $(X^{(\alpha)}_t)_{t\in[0,1]}$ takes the form
\[
R^{(\alpha)}(t,s)
=\begin{cases}
-\sqrt{(1-s)(1-t)}\,\ln\big(1-(s\wedge t)\big)
&\text{if }s,t\in[0,1),\\[2pt]
0&\text{if }s\in\{0,1\}\text{\ or }t\in\{0,1\}.
\end{cases}
\]
Performing the same steps as above, we obtain the differential equations
 \eqref{DE_e1}, \eqref{DE_e2} and \eqref{RBdiffeq} with $\alpha=1/2$.

According to Watson \cite[pp. 76, 83, 95]{Wat}, the functions $\sqrt{x}J_\nu(x)$ and $\sqrt{x}Y_\nu(x),$
 $x\in(0,1/\sqrt\lambda),$
 form a fundamental system of solutions for (\ref{RBdifeq}), and the general solution is of the form
\[
 f(x)=a_\nu\sqrt{x}\,J_\nu(x)
       +b_\nu\sqrt{x}\,Y_\nu(x),
\quad x\in(0,1/\sqrt\lambda),
\]
where $a_\nu$ and $b_\nu$ are arbitrary periodic functions of $\nu$ with period $1.$

Let us assume temporarily that $\alpha>1,$ i.e., $\nu>1/2.$ Then the boundary condition $ f(0)=0$ and
$\lim_{x\downarrow0}(\sqrt{x}\,J_\nu(x))=0$ (see part (i) of Proposition \ref{PRO:Bessel}) imply that
 $b_\nu\lim_{x\downarrow0}(\sqrt{x}\,Y_\nu(x))=0.$
However, since $\nu>1/2,$ by \eqref{Y_nu_kicsi_x}, we get
 \begin{align*}
  \lim_{x\downarrow 0}(\sqrt{x}\,Y_\nu(x))
    & = \lim_{x\downarrow 0}(x^{\frac{1}{2}-\nu} x^\nu Y_\nu(x))
      = \frac{-2^\nu\Gamma(\nu)}{\pi}
        \lim_{x\downarrow 0} x^{\frac{1}{2}-\nu}
      =-\infty.
 \end{align*}
Therefore $b_\nu$ must be zero and because of the periodicity,  $b_\nu=0$ even for $\nu\in(-1/2,1/2],$ i.e.,
 for $\alpha\in(0,1].$ Hence we can drop the assumption  $\alpha>1$ and obtain that $b_\nu=0$ for all
 $\nu>-1/2,$ i.e., for all $\alpha>0.$

So, the other boundary condition  $ f(1/\!\sqrt\lambda\,)=0$ implies
 $a_\nu J_\nu(1/\!\sqrt\lambda\,)=0.$
If $a_\nu=0,$ then $f(x)=0,$ $x\in(0,1/\sqrt{\lambda}),$ which yields that $e(t)=0,$ $t\in(0,1).$
However, since $e$ has to be an eigenfunction corresponding to the eigenvalue $\lambda,$
 it can not be identically zero and hence we have $a_\nu\ne0.$
This yields that $J_\nu(1/\!\sqrt\lambda\,)=0.$
Hence the eigenvalues $\lambda^{(\alpha)}_k=1/(z_k^{(\nu)})^{2},\ k\in\mathbb{N},$
 are the reciprocals of the squares of the (positive) zeros of $J_\nu.$
The corresponding normed eigenfunctions are
\[
 e^{(\alpha)}_k(t)
 =z^{(\nu)}_k
 \left(\int_0^{z^{(\nu)}_k}\!\!\!s\,J_\nu^2(s)\,\dd s\right)^{\!\!-\frac12}
 \!\sqrt{1-t}\,J_\nu\big(z_k^{(\nu)}(1-t)\big),\quad t\in[0,1],
\]
and using (\ref{5tulajd}) we obtain
\[
 e^{(\alpha)}_k(t)
 =\sqrt{2(1-t)}\,
 \frac{J_\nu\big(z_k^{(\nu)}(1-t)\big)}{\big\vert J_{\nu+1}\big(z_k^{(\nu)}\big)\big\vert}
 =\sqrt{2(1-t)}\,
 \frac{J_\nu\big(z_k^{(\nu)}(1-t)\big)}{\big\vert J_{\nu-1}\big(z_k^{(\nu)}\big)\big\vert},
 \quad t\in[0,1],
\]
 where we take the continuous extension of $e^{(\alpha)}_k$ at $t=1$ for $0<\alpha<1/2$, i.e.,
 \[
   e^{(\alpha)}_k(1)
      = \frac1{\big\vert J_{\nu+1}\big(z_k^{(\nu)}\big)\big\vert}
        \lim_{t\uparrow 1}\big(\sqrt{2(1-t)} J_\nu(z_k^{(\nu)}(1-t))\big)
      =0
 \]
 for $0<\alpha<1/2$ (see part (i) of Proposition \ref{PRO:Bessel}).
We also note that $J_{\nu+1}(z^{(\nu)}_k)\ne 0$ and $J_{\nu-1}(z^{(\nu)}_k)\ne 0,$
 by part (iii) of  Proposition \ref{PRO:Bessel}.

Now we return to the interval $[0,T]$.
Using (\ref{covfv}) (and recalling that we have already extended the covariance function
 $R^{(\alpha)}$ continuously onto $[0,T]^2$\!) we obtain the scaling law
 \begin{equation}
 \label{scalinglaw}
   R^{(\alpha,T)}(s,t)=TR^{(\alpha,1)}\left(\frac{s}{T},\frac{t}{T}\right),\quad s,t\in[0,T],
 \end{equation}
 where $R^{(\alpha,T)},$ resp. $R^{(\alpha,1)}$ denotes the covariance function on $[0,T]^2,$
 resp. on $[0,1]^2$ of the \aWb\ on the time interval $[0,T],$ resp. on $[0,1].$
Let $\lambda^{(\alpha,T)},$ resp. $\lambda^{(\alpha,1)}$ be a non-zero (and hence positive)
 eigenvalue of the integral operator $A_{R^{(\alpha,T)}},$ resp. $A_{R^{(\alpha,1)}}.$
Using (\ref{scalinglaw}) and the eigenvalue equation for $R^{(\alpha,T)},$
 one can write the eigenvalue equation for $R^{(\alpha,1)}$ in the following form
\begin{multline*}
\int_0^1\!\!R^{(\alpha,1)}(t,s) e^{(\alpha,T)}(Ts)\,\dd s
=\frac1T\int_0^1\!\! R^{(\alpha,T)}(Tt,Ts)e^{(\alpha,T)}(Ts)\,\dd s\\
=\frac1{T^2}\int_0^T\!\!R^{(\alpha,T)}(Tt,s) e^{(\alpha,T)}(s)\,\dd s
=\frac{\lambda^{(\alpha,T)}}{T^2}\,e^{(\alpha,T)}(Tt),\quad t\in[0,1],
\end{multline*}
 where the superscripts in parentheses denote the same correspondence as before.
 Similarly one can check that
 \begin{align*}
   \int_0^T\!\!R^{(\alpha,T)}(t,s) e^{(\alpha,1)}\left(\frac{s}{T}\right)\,\dd s
      = T^2 \lambda^{(\alpha,1)}e^{(\alpha,1)}\left(\frac{t}{T}\right),\quad t\in[0,T].
 \end{align*}
Then using part (iii) of Proposition \ref{PRO:Bessel} we have
 \[
  \lambda^{(\alpha,T)}_k = T^2\lambda^{(\alpha,1)}_k,\qquad k\in\NN,
 \]
 and by norming the eigenfunctions we obtain
 \[
 e^{(\alpha,T)}_k(t)=\frac{1}{\sqrt{T}}\,e^{(\alpha,1)}_k\left(\frac{t}{T}\right),\quad t\in[0,T],\;k\in\NN.
 \]
Indeed,
 \begin{align*}
    & \int_0^1 \left( e^{(\alpha,T)}_k(Tt)\right)^2\,\dd t
        = \frac{1}{T}\int_0^T \left( e^{(\alpha,T)}_k(s)\right)^2\,\dd s
        =\frac{1}{T},\qquad k\in\NN,\\
    & \int_0^1 \left( e^{(\alpha,1)}_k(t)\right)^2\,\dd t = 1,\qquad k\in\NN,
 \end{align*}
 and then taking into account also that normed eigenfunctions are unique only up to sign
 we get $\sqrt{T} e^{(\alpha,T)}_k(Tt) = e^{(\alpha,1)}_k(t),$ $t\in[0,1],$ $k\in\NN.$
Hence the statement of the theorem follows.
\proofend

In the next remark we study the question whether $0$ is an eigenvalue of the integral operator
 $A_{R^{(\alpha)}}$ or not.

\begin{Rem}
We note that $0$ is not an eigenvalue of the integral operator $A_{R^{(\alpha)}}.$
Indeed, on the contrary let us suppose that $0$ is an eigenvalue of $A_{R^{(\alpha)}}.$
We may assume without loss of generality that $T=1$ (see the end of the proof of Theorem
 \ref{KLthm}).
Then there exists a function $e:[0,1]\to \RR$ which is not $0$ almost everywhere
 and
 \[
    \int_0^1 R^{(\alpha)}(t,s) e^{(\alpha)}(s)\,\dd s =0,
     \qquad t\in[0,1].
 \]
First let us suppose that $\alpha\neq1/2$. By \eqref{seged1}, we have
 \begin{align*}
    \int_0^1 (1-s)^\alpha e(s)\,\dd s
       - \int_0^t (1-s)^{1-\alpha} e(s)\,\dd s
       - (1-t)^{1-2\alpha}\int_t^1 (1-s)^{\alpha} e(s)\,\dd s
     = 0
 \end{align*}
 for all $t\in(0,1),$ and differentiating with respect to $t,$
 \begin{align*}
   - (1-t)^{1-\alpha} e(t)
     + (1-2\alpha)(1-t)^{-2\alpha}\int_t^1 (1-s)^{\alpha} e(s)\,\dd s
     + (1-t)^{1-2\alpha} (1-t)^\alpha e(t)
    =0,
 \end{align*}
 or equivalently
 \[
   \int_t^1 (1-s)^{\alpha} e(s)\,\dd s = 0, \qquad t\in(0,1).
 \]
Differentiating again with respect to $t$ one can derive $e(t)=0,$ $t\in(0,1),$ which leads us to a
 contradiction. The case $\alpha=1/2$ can be handled in a similar way.
\proofend
\end{Rem}

In the next remark we will study the convergence of the coefficients of the random variables in the terms
 on the right-hand side of the KL representation \eqref{KLrepr} as \ $\alpha\downarrow 0$.

\begin{Rem}
\label{KL0}
First we recall that in the case of $\alpha=0$ the process $(X^{(0)}_t)_{t\in[0,T)}$ is the standard
 Wiener process on $[0,T]$.
If $\alpha\downarrow 0,$  then the left-hand side of the KL representation (\ref{KLrepr})
 converges in $L^2(\Omega,\cA,P)$ to the standard Wiener process $(B_t)_{t\in[0,T]},$ uniformly in $t\in[0,S]$ on every interval $[0,S]\subset[0,T).$
Indeed, for all $t\in[0,S],$
 \begin{align*}
   \EE(X^{(\alpha)}_t - B_t)^2
      & =  \EE\left(\int_0^t\!\left( \left(\frac{T-t}{T-s}\right)^\alpha -1 \right)\dd B_s \right)^{\!2}\\[2mm]
      & = \int_0^t\!\left( \left(\frac{T-t}{T-s}\right)^\alpha -1 \right)^{\!2}\dd s
          \leq S\left(\left(\frac{T-S}T\right)^{\!\alpha}-1\right)^{\!2}
         \to 0 \quad \text{as $\alpha\downarrow 0,$}
 \end{align*}
 where the last inequality follows by $(T-S)/T\leq (T-t)/(T-s)\leq 1,$ $0\leq s\leq t\leq S<T.$
Hence
 \[
   \sup_{t\in[0,S]} \EE(X^{(\alpha)}_t - B_t)^2 \to 0 \quad \text{as $\alpha\downarrow 0.$}
 \]
Hereafter we show that the coefficients of the random variables in the terms on the right-hand
 side of (\ref{KLrepr}) also converge uniformly in $t\in[0,S]$ to those of the corresponding terms
 of the KL expansion of $(B_t)_{t\in[0,T]},$ based on  $[0,T].$
For the latter KL representation see, e.g., Papoulis \cite[Example 12.10]{Papo}
 (which unfortunately contains a misprint).
Indeed, exploiting the fact that the eigenfunctions are unique only up to sign,
the KL expansion of $(B_t)_{t\in[0,T]},$ based on  $[0,T],$ can be written in the form
\begin{equation}
\label{WienerKL}
B_t=\sum_{k=1}^\infty\eta_k\,(-1)^{k-1}\sqrt{2T}\,\frac{\sin\big((k-1/2)\pi t/T\big)}
    {(k-1/2)\pi}
\end{equation}
 for all $t\in[0,T],$ where $\eta_k,\ k\in\mathbb{N},$ are independent, standard normally
 distributed random variables.
Moreover, using Theorem \ref{KLthm}, parts (ii) and (vii) of Proposition \ref{PRO:Bessel}
 and that therefore $z_k^{(-1/2)}=(k-1/2)\pi,\ k\in\mathbb{N},$  we obtain
\begin{align*}
\lim_{\alpha\downarrow0}\left(\sqrt{\lambda^{(\alpha)}_k} e^{(\alpha)}_k(t)\right)
&=\lim_{\alpha\downarrow0}
\left(\frac{T}{z_k^{(\nu)}}\,\sqrt{\frac{2}T\,\left(1-\frac{t}T\right)}
\,\frac{J_\nu\big(z_k^{(\nu)}(1-t/T\big)}{\big\vert J_{\nu+1}\big(z_k^{(\nu)}\big)\big\vert}\right)\\[5pt]
&=\frac{T}{z_k^{(-1/2)}}\,\sqrt{\frac{2}T\,\left(1-\frac{t}T\right)}
\,\frac{J_{-1/2}\big(z_k^{(-1/2)}(1-t/T\big)}{\big\vert J_{1/2}\big(z_k^{(-1/2)}\big)\big\vert}\\[5pt]
&=\sqrt{2T}\
    \frac{\cos\big((k-1/2)\pi(1-t/T)\big)}
    {(k-1/2)\pi}\\[5pt]
&=(-1)^{k-1}\sqrt{2T}\,\frac{\sin\big((k-1/2)\pi t/T\big)}
    {(k-1/2)\pi}
\end{align*}
for each $k\in\mathbb{N}.$
Further, the convergence is uniform in $t\in[0,S],$
because $\lim_{\alpha\downarrow 0}z_k^{(\nu)}=\lim_{\nu\downarrow -1/2}z_k^{(\nu)}=z_k^{(-1/2)}>0,$
 and the function
 \[
    \left[z_k^{(-1/2)}\left(1-\frac{S}{T}\right)-\varepsilon,z_k^{(-1/2)}+\varepsilon\right]
       \times \left[-\frac{1}{2},-\frac{1}{2}+\varepsilon\right]\ni (x,\nu)
       \mapsto J_\nu(x)
 \]
 is uniformly continuous (where $\varepsilon>0$ is sufficiently small),
 since, by part (ii) of Proposition \ref{PRO:Bessel},
 $(0,\infty)\times(-1,\infty)\ni(x,\nu)\mapsto J_\nu(x)$ is an analytic and hence
 continuous function,
 $\left[z_k^{(-1/2)}\left(1-\frac{S}{T}\right)-\varepsilon,z_k^{(-1/2)}+\varepsilon\right]
       \times \left[-\frac{1}{2},-\frac{1}{2}+\varepsilon\right]$ is a compact set
 and a continuous function on a compact set is uniformly continuous.
\proofend
\end{Rem}

In the next remark we consider the special case $\alpha=1$ in Theorem \ref{KLthm}.

\begin{Rem}\label{KL1}
If $\alpha=1,$ i.e., $\nu=1/2,$ then by Theorem \ref{KLthm},
 part (vii) of Proposition \ref{PRO:Bessel} and that therefore $z_k^{(1/2)}=k\pi,\ k\in\mathbb{N},$
we obtain that the eigenvalue--normed eigenfunction pairs are
 $\lambda^{(1)}_k=T^2/(k\pi)^2$ and
 \begin{align*}
  e^{(1)}_k(t)
   & =\sqrt{\frac2T\left(1-\frac{t}T\right)}\,\frac{J_{1/2}\big(k\pi(1-t/T)\big)}
           {\big\vert J_{-1/2}(k\pi)\big\vert}\\[7pt]
   & =\sqrt{\frac2T\left(1-\frac{t}T\right)}\,
   \frac{\sqrt{2T\big/\big(k\pi^2(T-t)\big)}\,\sin\big(k\pi(1-t/T)\big)}
        {\sqrt{2\big/\big(k\pi^2\big)}\,\big\vert\cos(k\pi)\big\vert}\\[4pt]
   & = (-1)^k\sqrt{\frac{2}{T}}\sin\left(\frac{k\pi t}{T}\right),
    \qquad t\in[0,T],\ k\in\mathbb{N}.
 \end{align*}
Further,
\[
 X_t^{(1)} = \sqrt{2T}\,\sum_{k=1}^\infty \xi_k (-1)^k
             \frac{\sin\big(k\pi t/T\big)}{k\pi}
           \distre
            \sqrt{2T}\,\sum_{k=1}^\infty \xi_k
             \frac{\sin\big(k\pi t/T\big)}{k\pi}
\]
 for all $t\in[0,T],$ where $\distre$ denotes equality in distribution.
As we expected, this is the KL expansion of the Wiener bridge, see, e.g., Deheuvels \cite[Remark 1.1]{Deh1}
 or Guti\'errez and Valderrama \cite[formula (10)]{GutVal}.
\proofend
\end{Rem}

In the next remark we present a space-time transformed Wiener process representation of the \aWb,
needed further on.
The idea comes from the similar representation of the Wiener bridge, see, e.g., Cs\"org\H o and
 R\'ev\'esz \cite[Proposition 1.4.2]{CsorRev}, and from Barczy and Pap \cite[proof of Lemma 3.1]{BarPap1}.
For historical fidelity we note that our representation (for the \aWb) is an analogue
 of formula (20) in Brennan and Schwartz \cite{BreSch}.

\begin{Rem}\label{KL2}
Let $( W_u)_{u\geq0}$ be a standard Wiener process,
\begin{equation}
\label{tauT}
 \tau_T^{(\alpha)}(t) := \int_0^t\frac{1}{(T-s)^{2\alpha}}\,\dd s
  = \frac{R^{(\alpha)}(t,t)}{(T-t)^{2\alpha}},\qquad t\in[0,T)
\end{equation}
and
\begin{equation} \label{sptime-tr}
  Z_t^{(\alpha)}:=(T-t)^\alpha  W_{\tau_T^{(\alpha)}(t)}, \quad t\in[0,T).
\end{equation}
Since $\tau_T^{(\alpha)}(0)=0$ and  $\tau_T^{(\alpha)}$ is strictly increasing and continuous,
 $(Z_t^{(\alpha)})_{t\in[0,T)}$ can be called a space-time transformed Wiener process.
One can see at once that this is a Gaussian process with almost surely continuous sample paths,
zero mean and
the covariance function (\ref{covfv}), therefore the process $(Z_t^{(\alpha)})_{t\in[0,T)}$
is a weak solution of the SDE \eqref{alpha_W_bridge}.
Since the SDE \eqref{alpha_W_bridge} has a strong solution which is pathwise unique, we get
 $(Z_t^{(\alpha)})_{t\in[0,T)}$ is  an \aWb \ (there exists some appropriate standard Wiener process
 for which \eqref{MArepr} holds).
 \proofend
\end{Rem}

In the following we deal with the weighted KL expansion of the \aWb.
The series expansion which we call the weighted KL expansion of a space-time transformed centered
 process with continuous covariance function was introduced by Guti\'errez and Valderrama \cite{GutVal}.
Let $S\in(0,T)$ and $\mu_T^{(\alpha)}$ be a (weight) measure defined on  (the Borel sets of) $[0,S]$
by the help of the  space-time transform in (\ref{sptime-tr}) as
 \[
    \dd \mu_T^{(\alpha)}(s):=(T-s)^{-2\alpha}\dd\tau_T^{(\alpha)}(s)=(T-s)^{-4\alpha}\dd s
 \]
 and let us denote by $L^2\big([0,S],\mu_T^{(\alpha)}\big)$ the Hilbert space of measurable functions on $[0,S],$
which are square integrable with respect to $\mu_T^{(\alpha)}.$
Furthermore, let
\begin{equation}
\label{Wienerf-KLrepr}
 W_u=\sum_{k=1}^\infty\sqrt{\kappa^{(\alpha)}_k}\,\xi_k d^{(\alpha)}_k(u),\quad u\in[0,\tau_T^{(\alpha)}(S)],
\end{equation}
 be the (unweighted) KL expansion of the standard Wiener process $( W_u)_{u\in[0,\tau_T^{(\alpha)}(S)]},$
 based on $[0,\tau_T^{(\alpha)}(S)],$ i.e.,
 $(\kappa^{(\alpha)}_k,d^{(\alpha)}_k),\ k\in\mathbb{N},$ are the eigenvalue--normed eigenfunction
 pairs of the integral operator associated to the covariance function of the standard Wiener process
 (for explicit formulae see \eqref{lambdak} and \eqref{lambdak2} later on) and
 $\xi_k,\ k\in\mathbb{N},$ are independent standard normally distributed random variables.
Finally, let
\begin{equation} \label{wghteigfunct}
  f^{(\alpha)}_k(t):=(T-t)^\alpha d_k^{(\alpha)}(\tau_T^{(\alpha)}(t)),\quad t\in[0,S],\,\ k\in\mathbb{N},
\end{equation}
i.e., we apply the same time change and rescaling to the normed eigenfunction $d_k^{(\alpha)}$ in order to
 define $f_k^{(\alpha)}$ what we apply to a standard Wiener process in order to get an $\alpha$-Wiener bridge,
 see \eqref{sptime-tr}.
Using (\ref{sptime-tr}) and Guti\'errez and Valderrama \cite[Theorem 1]{GutVal},
 the weighted KL expansion of the \aWb\ $(X_t^{(\alpha)})_{t\in[0,T]},$ based on $[0,S],$ with respect
 to the weight measure $\mu_T^{(\alpha)}$ is
\begin{equation} \label{wghtKLrepr}
 X_t^{(\alpha)}=\sum_{k=1}^\infty\sqrt{\kappa^{(\alpha)}_k}\,\xi_k f^{(\alpha)}_k(t),\quad t\in[0,S].
\end{equation}
It also follows that the properties of the weighted KL expansion (\ref{wghtKLrepr}) and
the weighted normed eigenfunctions (\ref{wghteigfunct}) are completely analogous to those of
(\ref{Wienerf-KLrepr}) and the unweighted normed eigenfunctions therein.
The difference is in the measure with respect to which we integrate. Namely,
 in the weighted case we integrate with respect to a Lebesgue--Stieltjes measure and there is the
 $L^2\big([0,S],\mu_T^{(\alpha)}\big)$ space in the background,
instead of the Lebesgue measure and the $L^2([0,\tau_T^{(\alpha)}(S)])$ space in the unweighted case.
So, the series in \eqref{wghtKLrepr} is convergent in $L^2(\Omega,\cA,P)$ uniformly in $t\in[0,S]$
 and
 \begin{align}\nonumber
   & \int_0^S R^{(\alpha)}(t,s) f^{(\alpha)}_k(s)\,\dd \mu_T^{(\alpha)}(s)
       = \kappa^{(\alpha)}_kf^{(\alpha)}_k{(t)}, \qquad t\in[0,S],\;k\in\NN,\\\label{seged7}
   & \int_0^S (f^{(\alpha)}_k(s))^2 \,\dd \mu_T^{(\alpha)}(s) = 1,
      \quad \int_0^S f^{(\alpha)}_k(s) f^{(\alpha)}_\ell(s) \,\dd \mu_T^{(\alpha)}(s) = 0,
      \quad k\ne\ell,\;k,\ell\in\NN.
 \end{align}

By Papoulis \cite[Example 12.10]{Papo} (which unfortunately contains a misprint), we have
\begin{align}\label{lambdak}
  \kappa^{(\alpha)}_k&=\bigg(\frac{\tau_T^{(\alpha)}(S)}{(k-1/2)\pi}\bigg)^2,
               \quad k\in\mathbb{N},\\[5pt]\label{lambdak2}
  d^{(\alpha)}_k(u)&=\sqrt{\frac2{\tau_T^{(\alpha)}(S)}}\,\sin\!\left(\!\bigg(k-\frac12\bigg)
          \pi\frac{u}{\tau_T^{(\alpha)}(S)}\right),
 \quad u\in[0,\tau_T^{(\alpha)}(S)],\,\ k\in\mathbb{N},
\end{align}
and then using (\ref{wghteigfunct}--\ref{wghtKLrepr}) we obtain the following theorem.

\begin{Thm}\label{wghtKLthm}
\hskip-2pt  In the weighted KL expansion \eqref{wghtKLrepr} of the $\alpha$\!-Wiener bridge
 the weighted eigenvalues $\kappa^{(\alpha)}_k,\ k\in\mathbb{N},$ are given by \eqref{lambdak}
 and the corresponding weighted normed eigenfunctions
\begin{equation}
\label{ekk}
 f^{(\alpha)}_k(t)
   =\sqrt{\frac2{\tau_T^{(\alpha)}(S)}}\,(T-t)^\alpha\sin\!\left(\!\bigg(k-\frac12\bigg)
           \pi\,\frac{\tau_T^{(\alpha)}(t)}{\tau_T^{(\alpha)}(S)}\right),
\quad t\in[0,S],\,\ k\in\mathbb{N}.
\end{equation}
Hence
 \[
   X^{(\alpha)}_t
      = \sum_{k=1}^\infty \sqrt{2\tau_T^{(\alpha)}(S)} (T-t)^\alpha
        \frac{\sin\left(\left(k-\frac{1}{2}\right)\pi\frac{\tau_T^{(\alpha)}(t)}{\tau_T^{(\alpha)}(S)}\right)}
             {(k-1/2)\pi}\xi_k,
          \qquad t\in[0,S].
 \]
\end{Thm}

In the next remark we will study the convergence of the coefficients of the random variables in the terms
 on the right-hand side of (\ref{wghtKLrepr}) as $\alpha\downarrow 0$.

\begin{Rem}\label{wghtKLrem1}
If $\alpha\downarrow 0,$ then the left-hand side of the weighted KL representation (\ref{wghtKLrepr})
 converges in $L^2(\Omega,\cA,P)$ to $B_t$ uniformly in $t\in[0,S],$ see the beginning of Remark \ref{KL0}.
Hereafter we show that the coefficients of the random variables in the terms on the right-hand side of
 (\ref{wghtKLrepr}) (given by the help of \eqref{tauT}, \eqref{lambdak} and \eqref{ekk})
 converge uniformly in $t\in[0,S]$ to the coefficients of the corresponding terms of
 the (unweighted) KL expansion of the standard Wiener process, based on $[0,S].$
Indeed, we have
\[
\lim_{\alpha\downarrow 0}\sqrt{2\tau_T^{(\alpha)}(S)}\,(T-t)^{\alpha}=\sqrt{2S},
\]
uniformly in $t\in[0,S]$ (the uniform convergence follows by mean value theorem), and
\[
\lim_{\alpha\downarrow0}\frac{\tau_T^{(\alpha)}(t)}{\tau_T^{(\alpha)}(S)}
=\frac{t}{S}\,,
\]
also uniformly in $t\in[0,S].$
Hence for all $k\in\NN,$
\begin{align*}
\lim_{\alpha\downarrow0}
\sqrt{\kappa^{(\alpha)}_k}\, f^{(\alpha)}_k(t)
&=\lim_{\alpha\downarrow0}
\sqrt{2\tau_T^{(\alpha)}(S)}\,(T-t)^{\alpha}
               \frac{\sin\big((k-1/2)\pi\tau_T^{(\alpha)}(t)\big/\tau_T^{(\alpha)}(S)\big)}{(k-1/2)\pi}\\[5pt]
&=\sqrt{2S}\,\frac{\sin\big((k-1/2)\pi t/S\big)}{(k-1/2)\pi},
\end{align*}
uniformly in $t\in[0,S].$
Indeed,
 \begin{align*}
    \left\vert
          \sin\left(\left(k-\frac{1}{2}\right)\pi\frac{\tau_T^{(\alpha)}(t)}{\tau_T^{(\alpha)}(S)}\right)
          - \sin\left(\left(k-\frac{1}{2}\right)\pi\frac{t}{S}\right)
    \right\vert
    \leq \left(k-\frac{1}{2}\right)\pi
         \left\vert\frac{\tau_T^{(\alpha)}(t)}{\tau_T^{(\alpha)}(S)} - \frac{t}{S} \right\vert.
 \end{align*}
\proofend
\end{Rem}

In the next remark we consider the special case $\alpha=1$ in Theorem \ref{wghtKLthm}.

\begin{Rem}\label{wghtKLrem2}
If $\alpha=1,$ then $\tau_T^{(\alpha)}$  (given by \eqref{tauT}) takes the form
 $\tau_T^{(1)}(t)=t/(T(T-t)),$ $t\in[0,T),$ so (\ref{lambdak}) takes the form
\[
 \kappa^{(1)}_k=\bigg(\frac{S}{T(T-S)}\bigg)^2 \frac1{(k-1/2)^2\pi^2},
    \quad k\in\mathbb{N},
\]
and (\ref{ekk}) becomes
\[
 f^{(1)}_k(t)=\sqrt{\frac{2T(T-S)}S}\,
          (T-t)\sin\left(\!\left(k-\frac12\right)\pi\,\frac{t(T-S)}{S(T-t)}\right),
\quad t\in[0,S],\,\ k\in\mathbb{N}.
\]
Particularly, for $T=1$ and $S=1/2$ we reobtain the weighted KL expansion
\[
X^{(1)}_t=
  \sqrt2\,(1-t)\sum_{k=1}^\infty\xi_k \frac{\sin\big((k-1/2)\pi t/(1-t)\big)}{(k-1/2)\pi},
  \quad t\in[0,1/2],
\]
 given by Guti\'errez and Valderrama \cite[formula (12)]{GutVal}.
\proofend
\end{Rem}

In the next remark we formulate a corollary of Theorem \ref{wghtKLthm} in the case of $\alpha=1/2$.

\begin{Rem}\label{wghtKLrem3}
For all \ $0<S<T$, \ we have
 \begin{align}\label{seged9}
      \int_0^S \frac{(X_u^{(1/2)})^2}{(T-u)^2}\,\dd u
              = \left(\ln\left(\frac{T}{T-S}\right)\right)^2
                \sum_{k=1}^\infty\frac1{(k-1/2)^2\pi^2}\xi_k^2,
 \end{align}
 where \ $\xi_k,$ $k\in\NN,$ \ are independent standard normally distributed random variables.
Indeed, by Theorem \ref{wghtKLthm} and the Parseval identity in $L^2([0,S],\mu_T^{(1/2)}),$
 we get
 \[
    \int_0^S \frac{(X_u^{(1/2)})^2}{(T-u)^2}\,\dd u
    = \sum_{k=1}^\infty \kappa^{(1/2)}_k\xi_k^2
         = \left(\ln\left(\frac{T}{T-S}\right)\right)^2
           \sum_{k=1}^\infty \frac{1}{(k-1/2)^2\pi^2}\xi_k^2,
 \]
 where the last equality follows by
 \[
   \tau^{(1/2)}_T(t) = \int_0^t \frac{1}{T-u}\,\dd u
                     = \ln\left(\frac{T}{T-t}\right), \qquad t\in[0,T).
 \]
\proofend
\end{Rem}

\section{Applications}\label{appl}

In this section we present some applications of the KL expansion \eqref{KLrepr} given in Theorem \ref{KLthm}.
First we calculate the Laplace transform of the $L^2[0,T]$\!-norm square of $(X_t^{(\alpha)})_{t\in[0,T]}.$

\begin{Pro}\label{appl1}
Let $T>0,$ $\alpha>0$ and $\nu:=\alpha-1/2.$
Then
 \begin{align}\label{seged4}
   &\EE\exp\left\{-c \int_0^T (X^{(\alpha)}_u)^2\,\dd u \right\}
       = \prod_{k=1}^\infty \frac{1}{\sqrt{1+2c T^2/(z^{(\nu)}_k)^2}}\,,
       \qquad c\geq 0,\\[2mm] \label{seged5}
   &\EE\exp\left\{-c \int_0^1 (X^{(0)}_u)^2\,\dd u \right\}
       = \frac{1}{\sqrt{\cosh(\sqrt{2c})}}, \qquad c\geq 0, \\[2mm] \label{seged6}
   & \EE\exp\left\{-c \int_0^1 (X^{(1)}_u)^2\,\dd u \right\}
       = \sqrt{\frac{\sqrt{2c}}{\sinh(\sqrt{2c})}}\,, \qquad c>0.
 \end{align}
Further, for all $0<S<T,$
 \begin{align}\label{seged8}
   & \EE\exp\left\{-c \int_0^S \frac{(X^{(1/2)}_u)^2}{(T-u)^2} \,\dd u \right\}
       = \frac{1}{\sqrt{\cosh\left(\sqrt{2c}\ln\left(\frac{T}{T-S}\right)\right)}},  \qquad c\geq 0.
 \end{align}
\end{Pro}

\noindent\textbf{Proof.}
By \eqref{KLrepr}, we have
 \begin{align*}
    (X^{(\alpha)}_t)^2 = \sum_{k=1}^\infty \sum_{\ell=1}^\infty
                           \sqrt{\lambda^{(\alpha)}_k \lambda^{(\alpha)}_\ell}
                           \,\xi_k\xi_\ell\, e^{(\alpha)}_k(t)e^{(\alpha)}_\ell(t),\quad t\in[0,T],
 \end{align*}
 and hence using that $\{e^{(\alpha)}_k,\, k\in\NN\}$ is an orthonormal system in $L^2[0,T]$ we get
 \begin{align}\label{seged3}
   \begin{split}
   \int_0^T  (X^{(\alpha)}_u)^2\,\dd u
      = \sum_{k=1}^\infty \sum_{\ell=1}^\infty
        \sqrt{\lambda^{(\alpha)}_k \lambda^{(\alpha)}_\ell} \,\xi_k\xi_\ell
          \int_0^T\!\! e^{(\alpha)}_k(u)e^{(\alpha)}_\ell(u)\,\dd u
      = \sum_{k=1}^\infty \lambda^{(\alpha)}_k\xi_k^2,
  \end{split}
 \end{align}
 which is nothing else but the Parseval identity in $L^2[0,T].$
Since $\xi_k,k\in\NN,$ are independent standard normally distributed random variables,
 for all $c\geq 0$ we get
 \begin{align*}
   \EE\exp\left\{-c \int_0^T (X^{(\alpha)}_u)^2\,\dd u \right\}
         = \prod_{k=1}^\infty \EE\left( \ee^{-c\lambda^{(\alpha)}_k\xi_k^2}\right)
         = \prod_{k=1}^\infty \frac{1}{\sqrt{1+2c\lambda^{(\alpha)}_k}},
 \end{align*}
 which yields \eqref{seged4}.

Now we turn to prove \eqref{seged5}.
First we show that for all $c\geq 0,$
 \begin{align}
   \lim_{\alpha\downarrow 0}
          \EE\exp\left\{-c \int_0^1 (X^{(\alpha)}_u)^2\,\dd u \right\}
       &=     \EE\exp\left\{-c \int_0^1 (X^{(0)}_u)^2\,\dd u \right\}, \label{seged45balo}\\
   \lim_{\alpha\downarrow 0}
       \prod_{k=1}^\infty \frac{1}{\sqrt{1+2c T^2/(z^{(\nu)}_k)^2}}
      &= \prod_{k=1}^\infty \frac{1}{\sqrt{1+2c T^2/(z^{(-1/2)}_k)^2}}.\label{seged45jobbo}
 \end{align}
First we check \eqref{seged45balo}.
By Cauchy--Schwarz's inequality  and Minkowski's inequality we get
\begin{align}
\EE\left\vert(X^{(\alpha)}_u)^2-(X^{(0)}_u)^2\right\vert
&=\EE\big(\left\vert X^{(\alpha)}_u-X^{(0)}_u\right\vert\left\vert
  X^{(\alpha)}_u+X^{(0)}_u\right\vert\big)\notag\\[4pt]
&\leq\left\Vert X^{(\alpha)}_u-X^{(0)}_u\right\Vert_{2}
  \big(\left\Vert X^{(\alpha)}_u\right\Vert_{2}+\left\Vert X^{(0)}_u\right\Vert_{2}\big),
  \; u\in[0,1),
\label{CStria}
\end{align}
where $\left\Vert\,\cdot\,\right\Vert_2$ denotes the $L^2(\Omega,\cA,\Prob)$--norm.
According to the beginning of Remark \ref{KL0},
\[
\lim_{\alpha\downarrow 0}\left\Vert X^{(\alpha)}_u-X^{(0)}_u\right\Vert_{2}
 =\lim_{\alpha\downarrow 0}\left\Vert X^{(\alpha)}_u-B_u\right\Vert_{2}
 =0
\]
uniformly in $u$ from any subinterval of $[0,1),$  where $(B_t)_{t\in[0,1]}$
 is a standard Wiener process.
Obviously, $\Vert X^{(0)}_u\Vert_{2}^2=\EE B_u^2=u\leq1$ if $u\in[0,1],$
and using \eqref{covfv}, we get for all $\alpha\in(0,1/4)$ and $u\in[0,1],$
\[
\left\Vert X^{(\alpha)}_u\right\Vert_{2}^2
=\EE(X^{(\alpha)}_u)^2
=\frac{(1-u)^{2\alpha}}{1-2\alpha}\big(1-(1-u)^{1-2\alpha}\big)
=\frac{(1-u)^{2\alpha}-(1-u)}{1-2\alpha}
\leq2.
\]
Hence, by (\ref{CStria}),
\[
\EE\left\vert(X^{(\alpha)}_u)^2-(X^{(0)}_u)^2\right\vert
   \leq3\left\Vert X^{(\alpha)}_u-X^{(0)}_u\right\Vert_{2}
   \leq 3(\Vert X^{(\alpha)}_u \Vert_{2}+ \Vert X^{(0)}_u\Vert_{2})
   \leq 9, \quad u\in[0,1].
\]
Now we can use Lebesgue's dominated convergence theorem to obtain
\begin{align*}
\lim_{\alpha\downarrow 0}
\EE\left\vert\int_0^1 (X^{(\alpha)}_u)^2\,\dd u - \int_0^1 (X^{(0)}_u)^2\,\dd u\right\vert
&\leq\lim_{\alpha\downarrow 0}
\int_0^1\EE\left\vert(X^{(\alpha)}_u)^2-(X^{(0)}_u)^2\right\vert\,\dd u\\[3pt]
&\leq3\lim_{\alpha\downarrow 0}\int_0^1 \left\Vert X^{(\alpha)}_u-X^{(0)}_u\right\Vert_{2}\,\dd u\\[2pt]
&=0\,.
\end{align*}
This is the $L^1\!$-convergence of the integrals  which yields
 their convergence in distribution and so, the convergence (\ref{seged45balo}) of their Laplace transforms.

To prove \eqref{seged45jobbo} we  rewrite the infinite product on the left hand side of
 \eqref{seged45jobbo} in the following form
 \[
     \prod_{k=1}^\infty \frac{1}{\sqrt{1+2c T^2/(z^{(\nu)}_k)^2}}
       = \exp\left\{-\frac{1}{2} \sum_{k=1}^\infty
          \ln\Bigg(1+2c\frac{T^2}{\big(z^{(\nu)}_k\big)^2}\Bigg) \right\}.
 \]
Each zero $z^{(\nu)}_k$ of $J_\nu$ is a strictly increasing and continuous function of $\nu\in(-1,\infty)$
 (see part (iii) of Proposition \ref{PRO:Bessel}), hence
\[
\lim_{\alpha\downarrow0}\ln\Bigg(1+2c\frac{T^2}{\big(z^{(\nu)}_k\big)^2}\Bigg)
 =\ln\Bigg(1+2c\frac{T^2}{\big(z^{(-1/2)}_k\big)^2}\Bigg),
 \quad k\in\NN,
\]
and for each $k\in\mathbb{N},$ the function
 \[
    (-1,\infty)\ni\nu \mapsto \ln\Bigg(1+2c\frac{T^2}{\big(z^{(\nu)}_k\big)^2}\Bigg)
 \]
 is decreasing, therefore the monotone convergence theorem can be applied to change the order of the
limit and the infinite sum, and after all, the order of the limit and the infinite product
on  the left-hand side of (\ref{seged45jobbo}).

The proof of \eqref{seged5} can be finished using \eqref{seged4} and Abramowitz and Stegun
 \cite[formula 4.5.69]{AbrSte}.  Indeed, we get
  \begin{align*}
    \EE\exp\left\{-c \int_0^1 (X^{(0)}_u)^2\,\dd u \right\}
           = \prod_{k=1}^\infty \frac{1}{\sqrt{1+2c\big/\big((k-1/2)^2\pi^2\big)}}
           = \frac{1}{\sqrt{\cosh(\sqrt{2c})}}\,.
  \end{align*}
This is a well-known formula due to L\'evy \cite{Lev}, see also Liptser and Shiryaev
 \cite[formula (7.147)]{LipShiI}.

Now we turn to prove \eqref{seged6}.
If $T=1$ and $\alpha=1$ (i.e., in the case of a Wiener bridge from $0$ to $0$ on the time interval $[0,1]$),
 by \eqref{seged4}, part (vii) of Proposition  \ref{PRO:Bessel} and
 Abramowitz and Stegun \cite[formula 4.5.68]{AbrSte},
 for all $c>0$,
  \begin{align*}
    \EE\exp\left\{-c \int_0^1 (X^{(1)}_u)^2\,\dd u \right\}
       = \prod_{k=1}^\infty \frac{1}{\sqrt{1+2c/(k^2\pi^2)}}
       = \sqrt{\frac{\sqrt{2c}}{\sinh(\sqrt{2c})}}\,.
  \end{align*}
This is also a well-known formula, see, e.g., Borodin and Salminen \cite[formula 1.9.7]{BorSal}.

Now we turn to prove \eqref{seged8}.
By  \eqref{seged9} we get for all $c\geq 0$
 \begin{align*}
   \EE\exp\left\{-c \int_0^S \frac{(X^{(1/2)}_u)^2}{(T-u)^2} \,\dd u \right\}
    & = \EE \prod_{k=1}^\infty\left(\exp\left\{ -c \left(\ln\left(\frac{T}{T-S}\right)\right)^2
                                    \frac{\xi_k^2}{(k-1/2)^2\pi^2}
                              \right\}\right)\\
    & = \prod_{k=1}^\infty
       \frac{1}{\sqrt{1+2c\frac{\left(\ln\left(\frac{T}{T-S}\right)\right)^2}{(k-1/2)^2\pi^2}}}
     =  \frac{1}{\sqrt{\cosh\left(\sqrt{2c}\ln\left(\frac{T}{T-S}\right)\right)}},
 \end{align*}
 where the last equality follows again by Abramowitz and Stegun \cite[formula 4.5.69]{AbrSte}.

We note that \eqref{seged8} yields that
 \[
    \PP\left(\lim_{S\uparrow T} \int_0^S \frac{(X^{(1/2)}_u)^2}{(T-u)^2} \,\dd u =\infty \right)=1,
 \]
 which was also obtained by Barczy and Pap \cite[(4.6)]{BarPap1} for general \aWb s.
Indeed,
 \[
   \lim_{S\uparrow T}\EE\exp\left\{-\int_0^S \frac{(X^{(1/2)}_u)^2}{(T-u)^2} \,\dd u \right\}
                     =0,
 \]
 and hence, by the Markov inequality,
 we get $\int_0^S \frac{(X^{(1/2)}_u)^2}{(T-u)^2} \,\dd u$ converges in probability
 to infinity as $S\uparrow T$.
Using Riesz's theorem and that the (random) function
 \[
     [0,T]\ni S \mapsto \int_0^S \frac{(X^{(1/2)}_u)^2}{(T-u)^2} \,\dd u
 \]
 is monotone increasing with probability one,  we get the desired property.
\proofend

We remark that a corresponding version of \eqref{seged8} for general \aWb s can be proved
 by a different technique, see Barczy and Pap \cite[Theorem 4.1]{BarPap2}.

Next we give a simple probabilistic proof for the sum of the square of the reciprocals of
 the positive zeros of $J_\nu$ with $\nu>-1/2.$
For $\nu>-1$ this is a well-known result due to Rayleigh, see, e.g., Watson \cite[Section 15.51, p. 502]{Wat}.
We note that Yor \cite[(11.47)--(11.49)]{Yor} and Deheuvels and Martynov \cite[Corollary 1.3]{DehMar}
 also gave probabilistic proofs of Rayleigh's results; we show that the proof of
 Deheuvels and Martynov can be carried through starting from the Karhunen--Lo\`eve expansion of the
\aWb\ as well.

\begin{Cor}\label{COR:Bessel}
Let $\alpha>0,$ $\nu:=\alpha-1/2$ and $z^{(\nu)}_k,\, k\in\NN,$ be the positive zeros of $J_\nu.$
Then
 \[
   \sum_{k=1}^\infty \frac{1}{(z^{(\nu)}_k)^2}
     = \frac{1}{4(\nu+1)}.
 \]
\end{Cor}

\noindent\textbf{Proof.}
Taking expectation of \eqref{seged3} we get
 \[
   \int_0^T  \EE(X^{(\alpha)}_u)^2\,\dd u
      = \sum_{k=1}^\infty \lambda^{(\alpha)}_k
      = T^2 \sum_{k=1}^\infty \frac{1}{(z^{(\nu)}_k)^2}.
 \]
By \eqref{covfv}, if $\alpha\ne 1/2$, then
 \begin{align*}
   \int_0^T  \EE(X^{(\alpha)}_u)^2\,\dd u
    & = \int_0^T \left(\frac{T^{1-2\alpha}}{1-2\alpha}(T-u)^{2\alpha} - \frac{T-u}{1-2\alpha}\right)
          \,\dd u
      = \frac{T^2}{1-4\alpha^2} - \frac{T^2}{2(1-2\alpha)}\\
    & = \frac{T^2}{2(1+2\alpha)},
 \end{align*}
 and if $\alpha=1/2$, then
 \begin{align*}
   \int_0^T  \EE(X^{(\alpha)}_u)^2\,\dd u
    & = \int_0^T (T-u)\ln\left(\frac{T}{T-u}\right) \,\dd u \\[1mm]
    &  =\ln(T)\int_0^T (T-u)\,\dd u - \int_0^T (T-u)\ln(T-u)\,\dd u \\[1mm]
    &  = \frac{T^2\ln(T)}{2}
         - \left(\frac{T^2\ln(T)}{2} - \frac{1}{2}\int_0^T (T-u)\,\dd u\right)
      =\frac{T^2}{4}.
 \end{align*}
Hence
 \[
   \sum_{k=1}^\infty \frac{1}{(z^{(\alpha-1/2)}_k)^2}
     = \frac{1}{2(1+2\alpha)}, \qquad \forall \; \alpha>0,
 \]
 which implies the statement.
\proofend

A consequence of (\ref{seged4}) is the following form of the survival function
 (complementary distribution function) of the $L^2[0,T]$\!-norm square of the \aWb.

\begin{Pro}\label{appl2} Let $\alpha>0$ and $\nu:=\alpha-1/2.$ Then
\label{compldf}
\begin{equation}
\label{compldf-forul}
   \P\left(\int_0^T\! (X^{(\alpha)}_t)^2\,\dd t > x\right)
   = \frac{2^{1-\nu/2}}{\pi\sqrt{\Gamma(\nu+1)}}
      \sum_{k=1}^\infty(-1)^{k+1}\!\!
      \int_{z_{2k-1}^{(\nu)}}^{z_{2k}^{(\nu)}}\!\!
      u^{\nu/2-1} \frac{\ee^{-xu^2/(2T^2)}}{\sqrt{\big\vert J_\nu(u)\big\vert}}\;\dd u
\end{equation}
 for all $x>0.$
\end{Pro}

\noindent\textbf{Proof.}
By Corollary \ref{COR:Bessel}, $\sum_{k=1}^\infty \lambda^{(\alpha)}_k<\infty$ and then
 \eqref{lambdakk} and the Smirnov formula (see, e.g., Smirnov \cite[formula (97)]{Smi}
 and Martynov \cite[formula (4)]{Mar1}, \cite[formula (2)]{Mar2}) imply that
\begin{equation}
\label{compleof}
   \P\left(\int_0^T\! (X^{(\alpha)}_t)^2\,\dd t > x\right)
   = \frac1\pi\sum_{k=1}^\infty(-1)^{k+1}\!\!
   \int_{(z_{2k-1}^{(\nu)})^2/T^2}^{(z_{2k}^{(\nu)})^2/T^2}
   \frac{\ee^{-ux/2}}{u\,\sqrt{\vert F(u)\vert}}\,\dd u, \qquad x>0,
\end{equation}
 where
\begin{equation}
\label{Fredholm}
F(u):=\prod_{k=1}^\infty\left(1-\lambda_k^{(\alpha)}\,u\right)
     =\prod_{k=1}^\infty\left(1-\frac{T^2}{\big(z_k^{(\nu)}\big)^2}\,u\right),\quad u\geq0,
\end{equation}
is the so-called Fredholm determinant.
Expressing this product by the help of the Euler formula (\ref{Euler})
we obtain
\begin{equation}
\label{Fredholmalak}
F(u)=J_\nu(\sqrt{u}\,T)\,\frac{\Gamma(\nu+1)}{\big(\sqrt{u}\,T/2\big)^\nu}\,,
 \qquad u>0.
\end{equation}
With this replacement (\ref{compleof}) becomes
\begin{align*}
\P\left(\int_0^T\! (X^{(\alpha)}_t)^2\,\dd t > x\right)
&= \frac1\pi\sum_{k=1}^\infty(-1)^{k+1}\!\!
   \int_{(z_{2k-1}^{(\nu)})^2/T^2}^{(z_{2k}^{(\nu)})^2/T^2}
   \frac{\ee^{-ux/2}\big(\sqrt{u}\,T/2\big)^{\nu/2}}
        {u\,\sqrt{\Gamma(\nu+1) \big\vert J_\nu(\sqrt{u}\,T)\big\vert}}\;\dd u\\[5pt]
&= \frac{2^{1-\nu/2}}{\pi\sqrt{\Gamma(\nu+1)}}
      \sum_{k=1}^\infty(-1)^{k+1}\!\!
      \int_{z_{2k-1}^{(\nu)}}^{z_{2k}^{(\nu)}}\!\!
      u^{\nu/2-1} \frac{\ee^{-xu^2/(2T^2)}}{\sqrt{\big\vert J_\nu(u)\big\vert}}\;\dd u
\end{align*}
 for all $x>0.$
\proofend

In the next remark we check that the formula for the survival function of the
 $L^2[0,1]$\!-norm square of a standard Wiener process (see, e.g.,
 Deheuvels and Martynov \cite[formula (1.50)]{DehMar}) can be derived by taking the limit
 of \eqref{compldf-forul}  with $T=1$ as $\alpha\downarrow 0.$

\begin{Rem}\label{appl3}
According to Deheuvels and Martynov \cite[formula (1.50)]{DehMar}
the survival function of the $L^2[0,1]$\!-norm square of
 a standard Wiener process is
\begin{equation}
\label{compleofW}
\P\left(\int_0^1\! B_t^2\,\dd t > x\right)
=\frac2\pi
      \sum_{k=1}^\infty(-1)^{k+1}\!\!
      \int_{(2k-3/2)\pi}^{(2k-1/2)\pi}\!\!
      \frac{\ee^{-xu^2/2}}{u\,\sqrt{-\cos(u)}}\;\dd u
\end{equation}
 for all $x>0.$
The right-hand side of (\ref{compleofW}) is continuous in $x\in(0,\infty),$
 which can be derived using Lebesgue's dominated convergence theorem.
Indeed,
 \begin{align*}
\int_{(2k-3/2)\pi}^{(2k-1/2)\pi}\!\!\frac{\ee^{-xu^2/2}}{u\,\sqrt{-\cos(u)}}\;\dd u
&\leq\frac{\ee^{-x(2k-3/2)^2\pi^2/2}}{(2k-3/2)\pi}
    \int_{(2k-3/2)\pi}^{(2k-1/2)\pi}\!\!\frac1{\sqrt{-\cos(u)}}\;\dd u\\[8pt]
&=\frac{\ee^{-x(2k-3/2)^2\pi^2/2}}{(2k-3/2)\pi}
    \int_{\pi/2}^{3\pi/2}\!\!\frac1{\sqrt{-\cos(u)}}\;\dd u,
     \quad k\in\NN,
\end{align*}
and, by D'Alembert's criterion, for all $x>0,$
 \[
   \sum_{k=1}^\infty \frac{\ee^{-x(2k-3/2)^2\pi^2/2}}{(2k-3/2)\pi} <\infty.
 \]
Then the left-hand side of (\ref{compleofW}) is also continuous in $x\in(0,\infty).$
Using the continuity of probability and that the $L^2[0,1]$-norm square of a standard Wiener process
 takes the value zero with probability $0,$ we have that the left-hand side of (\ref{compleofW}) is continuous
 in $x=0$ too, and
 \[
    \lim_{x\downarrow 0} \PP\left( \int_0^1\! B_t^2\,\dd t > x\right)
        = \PP\left( \int_0^1\! B_t^2\,\dd t > 0\right)
        = 1.
 \]
Hence $\P\big(\int_0^1 B_t^2\dd t > x\big)$ is continuous  at every $x\in\mathbb{R}.$
Using that $\int_0^1 (X^{(\alpha)}_t)^2\dd t$ converges in distribution
 to $\int_0^1 B_t^2\dd t$ as $\alpha\downarrow 0$ (which was verified
 in the proof of \eqref{seged45balo}), we get
\[
 \lim_{\alpha\downarrow0}\PP\left(\int_0^1\! (X^{(\alpha)}_t)^2\,\dd t > x\right)
   =\P\left(\int_0^1\! B_t^2\,\dd t > x\right)
\]
for all $x\in\mathbb{R}.$
Therefore the right-hand side of (\ref{compldf-forul}) must also converge to the right-hand side of
 (\ref{compleofW}) for every $x>0,$  i.e., the survival function of the
 $L^2[0,1]$\!-norm square of a standard Wiener process is the limit of the survival function of the
 $L^2[0,1]$\!-norm square of the $\alpha$-Wiener bridge (on the time interval $[0,1]$) as $\alpha\downarrow 0.$
\proofend
\end{Rem}

In the next remark we consider the Proposition \ref{compldf} with the special choices
 $\alpha=1$ and $T=1$.

\begin{Rem}\label{appl4}
With the special choices $\alpha=1,$ i.e., $\nu=1/2$  and $T=1$ in
 Proposition \ref{compldf} we have for all $x>0,$
\begin{align*}
\P\left(\int_0^1\! \big(X^{(1)}_t\big)^2\,\dd t > x\right)
&=\frac{2^{1-1/4}}{\pi\sqrt{\Gamma(3/2)}}
      \sum_{k=1}^\infty(-1)^{k+1}\!\!
      \int_{z_{2k-1}^{(1/2)}}^{z_{2k}^{(1/2)}}\!\!
      u^{1/4-1} \frac{\ee^{-xu^2/2}}{\sqrt{\big\vert J_{1/2}(u)\big\vert}}\;\dd u\\[5pt]
&=\frac{2^{3/4}}{\pi\sqrt{\sqrt{\pi}/2}}
      \sum_{k=1}^\infty(-1)^{k+1}\!\!
      \int_{(2k-1)\pi}^{2k\pi}\!\!
      u^{-3/4} \frac{\ee^{-xu^2/2}}{\sqrt{-\sqrt{2}\,\sin(u)/\sqrt{\pi u}}}\;\dd u\\[5pt]
&=\frac2\pi
      \sum_{k=1}^\infty(-1)^{k+1}\!\!
      \int_{(2k-1)\pi}^{2k\pi}\!\!
      \frac{\ee^{-xu^2/2}}{\sqrt{-u\sin(u)}}\;\dd u\,,
\end{align*}
where we used part (vii) of Proposition \ref{PRO:Bessel}.
We reobtained the survival function
 of the $L^2[0,1]$\!-norm square of the Wiener bridge, see Deheuvels and Martynov
 \cite[formula (1.51)]{DehMar}.
\proofend
\end{Rem}

\begin{Rem}
We note that Deheuvels and Martynov \cite[formula (1.43)]{DehMar} gave an expression for the survival
 function of the $L^2[0,1]\!$-norm square of a particularly weighted time transformed Wiener bridge,
 namely $(t^{1/2-\nu}X^{(1)}_{t^{2\nu}})_{t\in[0,1]},$ where $(X_t^{(1)})_{t\in[0,1]}$ is a Wiener bridge
 on the time interval $[0,1]$ and $\nu>0.$
We can notice that the only difference between that formula and our formula \eqref{compldf-forul}
 is the denominator of the fraction in the argument of the exponential function, namely
 instead of $4\nu$ we have $2T^2.$
This means that the distribution of the $L^2[0,1]\!$-norm square of the above mentioned particularly
 weighted time transformed Wiener bridge is the same as  the distribution of the $L^2[0,T]\!$-norm
 square of an  appropriate \aWb.
Namely, in case of $\alpha>1/2,$ i.e., $\nu>0,$ the random variables
 \[
   \int_0^1 t^{1-2\nu}(X^{(1)}_{t^{2\nu}})^2\,\dd t
     = \int_0^1 t^{2(1-\alpha)}(X^{(1)}_{t^{2\alpha-1}})^2\,\dd t
  \qquad \text{and}\qquad
  \int_0^{\sqrt{2\alpha-1}}(X^{(\alpha)}_t)^2\,\dd t
 \]
 have the same distribution,  where $(X^{(\alpha)}_t)_{t\in[0,\sqrt{2\alpha-1}]}$ is an
 $\alpha$-Wiener bridge on the time interval $[0,\sqrt{2\alpha-1}]$.
 \proofend
\end{Rem}

Zolotarev \cite[formula (6)]{Zol} gives the tail behaviour of the distribution function of
 $\sum_{k=1}^\infty\lambda_k\xi_k^2,$ where $\xi_k,$ $k\in\NN,$ are independent standard normally
  distributed random variables and $(\lambda_k)_{k\in\NN}$ is a sequence of positive real numbers such that
  $\lambda_1>\lambda_2>\cdots>0$ and $\sum_{k=1}^\infty\lambda_k<\infty$
 (see also Hwang \cite[Theorem 1]{Hwa} or Deheuvels and Martynov \cite[Lemma 1.1 and Remark 1.2]{DehMar}).
This result can be directly applied together with Theorem \ref{KLthm} to obtain the following corollary
about the large deviation probabilities for the $L^2[0,T]$\!-norm square of the \aWb.

\begin{Cor}\label{appl5}
Let $\alpha>0,$ $\nu:=\alpha-1/2$ and $z^{(\nu)}_k,$ $k\in\NN,$ be the positive zeros of $J_\nu.$
Then
\begin{align}
\P\left(\int_0^T\! (X^{(\alpha)}_t)^2\,\dd t > x\right)
 &=\big(1+o(1)\big)\frac{2^{1-\nu/2}T
  \big(z_1^{(\nu)}\big)^{(\nu-3)/2}}{\sqrt{\pi\,\Gamma(\nu+1)J_{\nu+1}\big(z_1^{(\nu)}\big)}}
    \ x^{-1/2}\,\ee^{-(z_1^{(\nu)})^2 x/(2T^2)} \label{elso}\\[1mm]
   &=\big(1+o(1)\big)\sqrt{\frac2\pi}\,\frac{T}{z_1^{(\nu)}}
      \prod_{k=2}^\infty\left(1-\frac{(z_1^{(\nu)})^2}{(z_k^{(\nu)})^2}\right)^{\!\!-1/2}\!\!
         x^{-1/2}\,\ee^{-(z_1^{(\nu)})^2 x/(2T^2)} \label{masodik}
\end{align}
as $x\to\infty.$
\end{Cor}

\noindent\textbf{Proof.}
Theorem \ref{KLthm} and Zolotarev \cite[formula (6)]{Zol} (or see also Hwang \cite[Theorem 1]{Hwa} or
 Deheuvels and Martynov \cite[Lemma 1.1 and Remark 1.2]{DehMar}) yield that for all $x>0$,
\begin{equation}
\label{ebbbe}
\P\left(\int_0^T\! (X^{(\alpha)}_t)^2\,\dd t > x\right)
=\big(1+o(1)\big)\sqrt{\frac2\pi}\,\frac{T^2}{(z_1^{(\nu)})^2\sqrt{-F'\big((z_1^{(\nu)})^2/T^2\big)}}\,
\,x^{-1/2}\,\ee^{-(z_1^{(\nu)})^2 x/(2T^2)},
\end{equation}
 where $F$ is defined in \eqref{Fredholm}.
Now we check that
\begin{align}\label{Fredholmderivalt}
F'\bigg(\frac{\big(z_1^{(\nu)}\big)^2}{T^2}\bigg)
 =-\Gamma(\nu+1)2^{\nu-1}T^2(z_1^{(\nu)})^{-\nu-1}J_{\nu+1}\big(z_1^{(\nu)}\big).
\end{align}
By \eqref{Fredholmalak} we get for all $u>0,$
\[
 F'(u)
    =\Gamma(\nu+1)2^{\nu-1}T^{-\nu}u^{-\frac{\nu}{2}-1}
      \big(-\nu J_\nu(\sqrt{u} T) + T\sqrt{u}J_\nu'(\sqrt{u}T)
       \big),
\]
 and so,
\begin{align*}
F'\bigg(\frac{(z_1^{(\nu)})^2}{T^2}\bigg)
 &=\Gamma(\nu+1)2^{\nu-1}T^{-\nu}
   \bigg(\frac{(z_1^{(\nu)})^2}{T^2}\bigg)^{-\frac{\nu}{2}-1}
    \big(-\nu J_\nu(z_1^{(\nu)}) + z_1^{(\nu)}J_\nu'(z_1^{(\nu)})\big) \\
 &=\Gamma(\nu+1)2^{\nu-1}T^2(z_1^{(\nu)})^{-\nu-1}J_\nu\hskip-2pt'\big(z_1^{(\nu)}\big).
\end{align*}
By Watson \cite[p. 45, (3)--(4)]{Wat}, the derivative of $J_\nu$ can be
 expressed as $J_\nu\hskip-2pt'(x)=\nu J_\nu(x)/x-J_{\nu+1}(x).$
This yields \eqref{Fredholmderivalt} and hence \eqref{elso}.

Now we turn to prove \eqref{masodik}. For this it is enough to check that
 $F'((z_1^{(\nu)})^2/T^2)$ can be also written in the form
 \begin{align}\label{Fredholmderivalt2}
 F'\bigg(\frac{\big(z_1^{(\nu)}\big)^2}{T^2}\bigg)
   =-\frac{T^2}{(z_1^{(\nu)})^2}\,\prod_{k=2}^\infty\bigg(1-\frac{(z_1^{(\nu)})^2}{(z_k^{(\nu)})^2}\bigg).
\end{align}
By \eqref{Fredholm}, we have
 \begin{align*}
    F(u)=(1-u\lambda_1^{(\alpha)})\prod_{k=2}^\infty (1-u\lambda_k^{(\alpha)}),
    \qquad u>0,
 \end{align*}
 and hence
 \begin{align*}
    F'(u)=-\lambda_1^{(\alpha)}\prod_{k=2}^\infty (1-u\lambda_k^{(\alpha)})
           +(1-u\lambda_1^{(\alpha)})\frac{\dd}{\dd u}\left(\prod_{k=2}^\infty (1-u\lambda_k^{(\alpha)})\right),
    \qquad u>0.
 \end{align*}
This yields that
 \[
    F'\bigg(\frac{1}{\lambda_1^{(\alpha)}}\bigg)
      = -\lambda_1^{(\alpha)}\prod_{k=2}^\infty
          \left(1-\frac{\lambda_k^{(\alpha)}}{\lambda_1^{(\alpha)}}\right),
 \]
 which implies \eqref{Fredholmderivalt2}  using \eqref{lambdakk}.
\proofend

The next corollary describes the behaviour at zero of the distribution function
 (small deviation probabilities) of the $L^2[0,T]$\!-norm square of the \aWb.

\begin{Cor}\label{appl6}
Let $\alpha>0$ and $\nu:=\alpha-1/2.$
Then there exists some constant $c>0$ such that
\begin{equation}
\label{corollequ}
\P\left(\int_0^T\! (X^{(\alpha)}_t)^2\,\dd t <\varepsilon\right)
=\big(c+o(1)\big)\varepsilon^{1/4-\nu/2}\ee^{-T^2/(8\varepsilon)}
\end{equation}
as $\varepsilon\downarrow 0.$
\end{Cor}

\noindent\textbf{Proof.}
We imitate the proof of Deheuvels and Martynov \cite[Theorem 1.7]{DehMar},
which originates from the idea of Li \cite{Li1}.
By (\ref{seged3}) and (\ref{lambdakk}) we have
\begin{equation}
\label{Parsevalegy}
\int_0^T (X^{(\alpha)}_t)^2\,\dd t
=T^2\sum_{k=1}^\infty\frac{\xi_k^2}{(z_k^{(\nu)})^2}
\end{equation}
 with independent standard normally distributed random variables  $\xi_k,\ k\in\NN.$
We are going to replace the zeros $z_k^{(\nu)}$ on the right-hand side of (\ref{Parsevalegy})
 by the first terms ('leading terms') on the right-hand side of (\ref{gyokkonverzio}).
We use Theorem 2 in Li \cite{Li1} to prove that this replacement can be done.
Namely, given any two convergent series
 $\sum_{k=1}^\infty a_k$ and $\sum_{k=1}^\infty b_k$ with positive terms,
 the convergence condition $\sum_{k=1}^\infty\vert1-a_k/b_k\vert<\infty$ implies that
 \[
  \P\left(\sum_{k=1}^\infty a_k\xi_k^2\leq\varepsilon\right)
    =\big(1+o(1)\big)\sqrt{\prod_{k=1}^\infty\frac{b_k}{a_k}}\
     \P\!\left(\sum_{k=1}^\infty b_k\xi_k^2\leq\varepsilon\right)
 \]
 as $\varepsilon\downarrow 0.$
Observe that $\sum_{k=1}^\infty\vert 1-a_k/b_k\vert<\infty$ yields the converge
 of $\prod_{k=1}^\infty(b_k/a_k)$.
Indeed, by mean value theorem,
 \[
     \ln\frac{b_k}{a_k}
       = \ln\frac{b_k}{a_k} - \ln1
       = \frac{1}{\theta_k}
         \left( \frac{a_k}{b_k} - 1 \right),
         \qquad k\in\NN,
 \]
 with some $\theta_k$ from the interval with endpoints $a_k/b_k$ and $1.$
Since $\lim_{k\to\infty}a_k/b_k=1$, we have that the sequence $\theta_k,\,k\in\NN,$ is bounded
 and hence the series $\sum_{k=1}^\infty \ln\frac{b_k}{a_k}$ is absolute convergent.
In what follows let $a_k:=\big(z_k^{(\nu)}\big)^{\!-2},$ $b_k:=\big(k+(\nu-1/2)/2\big)^{\!-2}\pi^{-2},$
 $k\in\mathbb{N}.$
Then by the help of (\ref{gyokkonverzio}) and using also that $z^{(\nu)}_1\leq z^{(\nu)}_k,\,k\in\NN,$
 we obtain
\[
\sum_{k=1}^\infty\left\vert1-\frac{a_k}{b_k}\right\vert
 =\sum_{k=1}^\infty\left\vert1-\bigg(\frac{k+(\nu-1/2)/2)\pi}{z_k^{(\nu)}}\bigg)^2\right\vert
 =\sum_{k=1}^\infty \frac{c_1}{k^2}
 <\infty
\]
 with some constant $c_1>0,$ which yields that
\begin{align}
 &\P\left(\sum_{k=1}^\infty \frac{\xi_k^2}{(z_k^{(\nu)})^2}\leq\frac\varepsilon{T^2}\right)\notag\\[5pt]
 &\qquad=\big(1+o(1)\big)\sqrt{\prod_{k=1}^\infty\bigg(\frac{k+(\nu-1/2)/2)\pi}{z_k^{(\nu)}}\bigg)^{-2}}\
  \P\!\left(\sum_{k=1}^\infty\frac{\xi_k^2}{\big(k+(\nu-1/2)/2\big)^2}\leq\frac{\varepsilon\pi^2}{T^2}\right)\notag\\[5pt]
 &\qquad=\big(c_2+o(1)\,\big)\P\!\left(\sum_{k=1}^\infty\frac{\xi_k^2}{\big(k+(\nu-1/2)/2\big)^2}
  \leq\frac{\varepsilon\pi^2}{T^2}\right)\label{ezkell}
\end{align}
 as $\varepsilon\downarrow 0,$ with a suitable constant $c_2>0.$
 By Li \cite[formula (3.4), p. 14]{Li1}, for an arbitrary $d>-1$
it holds that
\begin{equation}
\label{kiseltegyenl}
 \P\left(\sum_{k=1}^\infty \frac{\xi_k^2}{(k+d)^2}\leq\varepsilon\right)
  =\big(c_3+o(1)\big)
 \varepsilon^{-d}\ee^{-\pi^2/(8\varepsilon)}
\end{equation}
 as $\varepsilon\downarrow 0,$ with some constant $c_3>0.$
Combining (\ref{Parsevalegy}), (\ref{ezkell}) and (\ref{kiseltegyenl}) we obtain (\ref{corollequ}).
\proofend

\begin{Rem}
\label{Niki}
In case of $\alpha\geq 1/2,$ Corollary \ref{appl6} can be improved by which we mean that the constant $c$
 can be explicitly given.
Namely, by Nazarov \cite[Lemma 3.2]{Naz}, if $\xi_k,$ $k\in\NN,$ are independent standard normally distributed
 random variables, then for all $\nu\geq 0$,
 \begin{align*}
   \PP\left( \sum_{n=1}^\infty \frac{\xi_n^2}{(z_n^{(\nu)})^2}\leq \varepsilon^2 \right)
      \sim \sqrt{\frac{\sqrt{\pi}}{2^{\nu-1/2}\Gamma(1+\nu)}}
           \frac{\sqrt{2} \varepsilon_1^{1/2-\nu}}{\sqrt{\pi\cD_1}}
            \exp\left\{-\frac{\cD_1}{2\varepsilon_1^2}\right\},
 \end{align*}
 as $\varepsilon\to0,$ where
  \[
   \varepsilon_1=\varepsilon\sqrt{2\sin(\pi/2)}=\sqrt{2} \varepsilon
   \qquad \text{and} \qquad
   \cD_1 = \frac{1}{2\sin(\pi/2)}=1/2.
  \]
Hence
 \begin{align*}
   \PP\left( \sum_{n=1}^\infty \frac{\xi_n^2}{(z_n^{(\nu)})^2}\leq \varepsilon^2 \right)
      \sim \frac{2^{\frac{3}{2}-\nu} \pi^{-1/4}}{\sqrt{\Gamma(1+\nu)}}
           \varepsilon^{1/2-\nu} \ee^{-\frac{1}{8\varepsilon^2}},
         \qquad \text{as \ $\varepsilon\downarrow 0.$}
 \end{align*}
Then for all $T>0$ and $\varepsilon>0$ we have
 \begin{align*}
   \PP\left( T^2\sum_{n=1}^\infty \frac{\xi_n^2}{(z_n^{(\nu)})^2}\leq \varepsilon \right)
      \sim \frac{2^{\frac{3}{2}-\nu} \pi^{-1/4}}{\sqrt{\Gamma(1+\nu)}}
           \frac{\varepsilon^{\frac{1}{4}-\frac{\nu}{2}}}{T^{1/2-\nu}} \ee^{-\frac{T^2}{8\varepsilon}},
         \qquad \text{as \ $\varepsilon\downarrow 0.$}
 \end{align*}
By \eqref{Parsevalegy}, this yields that
  in case of $\alpha\geq 1/2$ the constant $c$ in Corollary \ref{appl6} takes the following form
  \[
     c = \frac{2^{\frac{3}{2}-\nu} \pi^{-1/4}}{\sqrt{\Gamma(1+\nu)}T^{1/2-\nu}}.
  \]
The reason for restricting ourselves to the case $\alpha\geq 1/2$ is that Lemma 3.2 in Nazarov \cite{Naz}
 is valid for $\nu\geq 0,$ while in Corollary \ref{appl6} we have $\nu=\alpha-1/2,$ $\alpha>0.$
\proofend
\end{Rem}

\section{Appendix}

In the next proposition we list some properties of the Bessel functions $J_\nu$ and $Y_\nu$ of
 the first and second kind, respectively (introduced in Section \ref{KLsect}).

\begin{Pro}\label{PRO:Bessel}
\renewcommand{\labelenumi}{{\rm(\roman{enumi})}}
\begin{enumerate}
\item For all $\nu\in\RR$, $J_\nu$ is continuous on $(0,\infty),$  and in case of $\nu\geq0$
      it can be continuously extended to $[0,\infty)$ by $J_\nu(0):=0$ if $\nu>0$ and by $J_\nu(0):=1$
      if $\nu=0.$
      However, in case of $-1/2<\nu<0$ we have $\lim_{x\downarrow 0} J_\nu(x)=\infty,$
      and $\lim_{x\downarrow 0}(\sqrt{x}J_\nu(x))=0$ holds for all $\nu>-1/2.$
      Further, for any $\nu>0,$
      \begin{align}\label{Y_nu_kicsi_x}
        Y_\nu(x) = (1+o(1))\frac{-\Gamma(\nu)}{\pi}\left(\frac{x}{2}\right)^{-\nu}
        \qquad \text{as $x\downarrow 0.$}
      \end{align}

\item By Watson \cite[p. 44]{Wat},
      $(0,\infty)\times(-1,\infty)\ni (x,\nu)\mapsto J_\nu(x)$ is an analytic function of both variables.

\item By the Bessel--Lommel theorem, see, e.g., Watson \cite[pp. 478, 482]{Wat},
     if $\nu>-1,$ then $J_\nu$ has infinitely many positive real zeros with multiplicities one.
     We denote them by $z^{(\nu)}_k\!,\ k\in\mathbb{N},$ where we assume
     \[
        0<z^{(\nu)}_1<z^{(\nu)}_2<\cdots<z^{(\nu)}_k<z^{(\nu)}_{k+1}<\cdots.
     \]
     By Watson \cite[p. 508]{Wat} or Korenev \cite[p. 96]{Kor},
     for fixed $k\in\NN,$ $z^{(\nu)}_k$ is a strictly increasing and continuous function of
     $\nu\in(-1,\infty).$

\item By Korenev \cite[p. 96]{Kor}, if $\nu>-1,$ then
\begin{equation}
\label{gyokkonverzio}
 z^{(\nu)}_k
   =\bigg(k+\frac12\Big(\nu-\frac12\Big)\bigg)\pi + O\!\left(\frac1k\right)
   \qquad \text{as \ $k\to\infty.$}
\end{equation}

\item If $\nu>-1,$ then $J_\nu$ can be written by the help of its
 zeros $z^{(\nu)}_k\!,\ k\in\mathbb{N},$ as
\begin{equation}
\label{Euler}
J_\nu(x)=\frac{\left(\frac{x}2\right)^\nu}{\Gamma(\nu+1)}\prod_{k=1}^\infty\left(1-\frac{x^2}{(z^{(\nu)}_k)^2}\right),
\quad x\in(0,\infty).
\end{equation}
This is the so called Euler formula, see Watson \cite[p. 498]{Wat}.

\item By Bowman \cite[p. 107]{Bow}, for $\nu>-1$ and for all zeros $z^{(\nu)}_k\!,\ k\in\mathbb{N},$
      of $J_\nu,$ it holds that
\begin{equation}
\label{5tulajd}
\int_0^{z^{(\nu)}_k}\!\!\!x J_\nu^2(x)\,\dd x
=\frac{\big(z^{(\nu)}_k\big)^2}2\, J_{\nu+1}^2\big(z^{(\nu)}_k\big)
=\frac{\big(z^{(\nu)}_k\big)^2}2\, J_{\nu-1}^2\big(z^{(\nu)}_k\big).
\end{equation}

\item The Bessel functions $J_{1/2}$ and $J_{-1/2}$ (of the first kind) can be expressed in closed forms
      by elementary functions:
  \[
    J_{1/2}(x)=\sqrt{\frac{2}{\pi x}}\sin(x) \qquad \text{and}\qquad
    J_{-1/2}(x)=\sqrt{\frac{2}{\pi x}}\cos(x),
    \qquad x\in(0,\infty),
   \]
 see, e.g., Watson \cite[pp.~54, 55]{Wat}.
\end{enumerate}
\end{Pro}

\noindent\textbf{Proof.}
We only have to check \eqref{Y_nu_kicsi_x}, since all the other statements
 are immediate consequences of the definitions or references were given for them.
If $\nu\in\NN$, then \eqref{Y_nu_kicsi_x} follows by Abramowitz and Stegun \cite[p. 360, 9.1.9]{AbrSte}.
If $\nu>0$ and $\nu\not\in\NN,$ then
 \begin{align*}
   \lim_{x\downarrow 0} (x^\nu Y_\nu(x))
     & =-\frac{1}{\sin(\pi \nu)} \lim_{x\downarrow 0} (x^\nu J_{-\nu}(x))\\
     & =-\frac{1}{\sin(\pi \nu)} \lim_{x\downarrow 0}
         \left( x^\nu
          \sum_{k=0}^\infty\frac{(-1)^k}{k!\,\Gamma(k-\nu+1)}\,\left(\frac{x}2\right)^{2k-\nu}
          \right)\\
     &= -\frac{2^\nu}{\sin(\pi \nu)} \lim_{x\downarrow 0}
        \left( \sum_{k=0}^\infty\frac{(-1)^k}{k!\,\Gamma(k-\nu+1)}\,\left(\frac{x}2\right)^{2k}\right)\\
     &= -\frac{2^\nu}{\sin(\pi \nu)\Gamma(-\nu+1)}
      = -\frac{2^\nu}{\sin(\pi \nu)(-\nu)\Gamma(-\nu)}
      = -\frac{2^\nu\Gamma(\nu)}{\pi},
 \end{align*}
 where the last equality follows by that
 \[
    \Gamma(z)\Gamma(-z) = -\frac{\pi}{z\sin(\pi z)},
    \qquad z>0.
 \]
\proofend

\section*{Acknowledgements}
 The first author has been supported by the Hungarian Scientific Research Fund under Grant No.~OTKA-T-079128.
 This work has been finished while M. Barczy was on a post-doctoral position at the Laboratoire de Probabilit\'es
 et Mod\`{e}les Al\'eatoires, University Pierre-et-Marie Curie, thanks to NKTH-OTKA-EU FP7 (Marie Curie action)
 co-funded 'MOBILITY' Grant No. MB08-A 81263.
We thank Yu. Nikitin for drawing our attention to the facts in Remark \ref{Niki}.

\end{document}